\pdfoutput=1
\RequirePackage{ifpdf}
\ifpdf 
\documentclass[pdftex]{sigma}
\else
\documentclass{sigma}
\fi

\usepackage{mathrsfs}
\usepackage{tikz}
\usetikzlibrary{matrix}

\numberwithin{equation}{section}

\newcommand{\rd}{{\rm d}}
\newcommand{\rG}{{\rm G}}

\newcommand{\bL}{{\bf L}}

\newcommand{\cA}{\mathcal{A}}

\newcommand{\cE}{\mathcal{E}}

\newcommand{\cH}{\mathcal{H}}

\newcommand{\cL}{\mathcal{L}}

\newcommand{\cO}{\mathcal{O}}

\newcommand{\cR}{\mathcal{R}}

\newcommand{\sA}{\mathscr{A}}
\newcommand{\sB}{\mathscr{B}}

\newcommand{\sE}{\mathscr{E}}
\newcommand{\sF}{\mathscr{F}}
\newcommand{\sG}{\mathscr{G}}

\newcommand{\sM}{\mathscr{M}}

\newcommand{\sO}{\mathscr{O}}

\newcommand{\sR}{\mathscr{R}}

\newcommand{\sU}{\mathscr{U}}

\newcommand{\fu}{{\mathfrak u}}
\newcommand{\fM}{{\mathfrak M}}
\newcommand{\N}{{\mathbb N}}
\newcommand{\Z}{{\mathbb Z}}
\newcommand{\Q}{{\mathbb Q}}
\newcommand{\R}{{\mathbb R}}
\newcommand{\C}{{\mathbb C}}
\newcommand{\fpu}{\mathfrak{pu}}

\newcommand{\SO}{{\rm SO}}

\newcommand{\SU}{{\rm SU}}
\newcommand{\GL}{\mathrm{GL}}

\newcommand{\SL}{\mathrm{SL}}
\newcommand{\U}{{\rm U}}
\newcommand{\PU}{{\mathbb P}{\rm U}}
\newcommand{\PGL}{{\mathbb P}\GL}

\renewcommand{\P}{{\mathbb P}}
\renewcommand{\epsilon}{\varepsilon}
\newcommand{\eps}{\epsilon}

\newcommand{\Irr}{\operatorname{Irr}}
\newcommand{\Coh}{\operatorname{Coh}}
\newcommand{\Aut}{\operatorname{Aut}}
\newcommand{\End}{\operatorname{End}}
\newcommand{\cEnd}{\cE{\rm nd}}
\newcommand{\sEnd}{\sE{\rm nd}}

\newcommand{\Hom}{\operatorname{Hom}}
\newcommand{\ch}{\operatorname{ch}}
\newcommand{\coker}{\operatorname{coker}}
\newcommand{\delbar}{\bar \del}
\newcommand{\del}{\partial}
\newcommand{\diag}{\operatorname{diag}}
\newcommand{\dvol}{\operatorname{dvol}}
\newcommand{\id}{\operatorname{id}}

\newcommand{\ind}{\mathop{\mathrm{index}}}
\newcommand{\into}{\hookrightarrow}
\newcommand{\iso}{\cong}

\newcommand{\RF}{{\rm RF}}

\newcommand{\tr}{\operatorname{tr}}
\newcommand{\rk}{\operatorname{rk}}

\newcommand{\Rep}{\operatorname{Rep}}
\newcommand{\Sym}{\operatorname{Sym}}
\newcommand{\Ext}{\operatorname{Ext}}
\newcommand{\sign}{\operatorname{sign}}
\newcommand{\an}{\mathrm{an}}

\def\wtilde{\widetilde}

\newtheorem{Theorem}{Theorem}[section]
\newtheorem{Corollary}[Theorem]{Corollary}
\newtheorem{Lemma}[Theorem]{Lemma}
\newtheorem{Proposition}[Theorem]{Proposition}
 { \theoremstyle{definition}
\newtheorem{Definition}[Theorem]{Definition}
\newtheorem{Remark}[Theorem]{Remark} }

\begin{document}
\allowdisplaybreaks

\newcommand{\arXivNumber}{1207.6938}

\renewcommand{\PaperNumber}{017}

\FirstPageHeading

\ShortArticleName{Rigid HYM Connections on Tautological Bundles over ALE Crepant Resolutions}

\ArticleName{Rigid HYM Connections on Tautological Bundles \\ over ALE Crepant Resolutions in Dimension Three}

\Author{Anda DEGERATU~$^\dag$ and Thomas WALPUSKI~$^\ddag$}

\AuthorNameForHeading{A.~Degeratu and T.~Walpuski}

\Address{$^\dag$~University of Freiburg, Mathematics Institute, Germany}
\EmailD{\href{mailto:anda.degeratu@math.uni-freiburg.de}{anda.degeratu@math.uni-freiburg.de}}
\URLaddressD{\url{https://home.mathematik.uni-freiburg.de/degeratu/}}

\Address{$^\ddag$~Massachusetts Institute of Technology, Department of Mathematics, USA}
\EmailD{\href{mailto:walpuski@mit.edu}{walpuski@mit.edu}}
\URLaddressD{\url{https://math.mit.edu/~walpuski/}}

\ArticleDates{Received June 02, 2015, in f\/inal form February 06, 2016; Published online February 15, 2016}

\Abstract{For $G$ a f\/inite subgroup of ${\rm SL}(3,{\mathbb C})$ acting freely on ${\mathbb C}^3{\setminus} \{0\}$ a crepant resolution of the Calabi--Yau orbifold ${\mathbb C}^3\!/G$ always exists and has the geometry of an ALE non-compact manifold.
  We show that the tautological bundles on these crepant resolutions admit rigid Hermitian--Yang--Mills connections.
  For this we use analytical information extracted from the derived category McKay correspondence
  of Bridgeland, King, and Reid~[\textit{J.~Amer. Math. Soc.} \textbf{14} (2001),
  535--554].
  As a consequence we rederive
  multiplicative cohomological identities on the crepant resolution using
  the Atiyah--Patodi--Singer index theorem.
  These results are dimension three analogues of Kronheimer and
  Nakajima's results~[\textit{Math. Ann.} \textbf{288} (1990), 263--307] in
  dimension two.}


\Classification{53C07; 14F05; 58J20}

\section{Introduction}
\label{sec:intro}

Let $G$ be a f\/inite subgroup of $\SL(n,\C)$ and let $\pi\colon X \to
\C^n\!/G$ be a crepant resolution.
Often $X$ comes equipped with a collection of so-called tautological vector
bundles $\sR_\rho$ indexed by the set $\Irr(G)$ of irreducible
representations of $G$.
When $n=2$,  Gonzalez-Sprinberg and Verdier~\cite{GonzalezSprinberg1983} discovered that these
vector bundles lie at the heart of the McKay
correspondence.
A geometrical interpretation of this result was given
by Kronheimer and Nakajima~\cite{Kronheimer1990} using asymptotically
locally Euclidean (ALE) hyperk\"ahler metrics on $X$ and inf\/initesimally
rigid Hermitian--Yang--Mills (HYM) connections on $\sR_\rho$.

In this paper we consider the case $n=3$ and $G$ a f\/inite subgroup of $\SL(3,\C)$ acting freely on $\C^3{\setminus} \{0\}$;
this condition is equivalent, in retrospect, with the existence of ALE crepant resolutions $\pi\colon X \to \C^3\!/G$.
From the classif\/ication of f\/inite subgroups of $\SL(3,\C)$ started by Blichfeldt~\cite{Blichfeldt1917} and completed by Yau and Yu~\cite{Yau-Stephan1993}, $G$ must be abelian\footnote{Such a group must also be isomorphic to $\Z_n$ with $n$ odd, see Remark~\ref{rem:fg}. However, for this work, it is only relevant that $G$ is abelian and that it acts freely on $\C^3{\setminus} \{0\}$.}.
From the work of Bridgeland, King, and Reid~\cite{Bridgeland2001} and Craw and
Ishii~\cite{Craw2004}, it is known that all projective crepant
resolutions of such~$\C^3\!/G$ can be constructed explicitly via GIT,
as moduli spaces $\sM_\theta$ of $G$-constellations with respect to a
generic rational stability parameter $\theta$.
Moreover, they are naturally equipped with a collection of tautological line bundles $\{ \sR_\rho\colon \rho \in \Irr(G) \}$.
Because of the relation between GIT and K\"ahler reduction, the
analytif\/ication\footnote{For
  more details on analytif\/ication, the passage from the algebraic to the analytic category, we recommend Neeman's book~\cite{Neeman2007}.} $M_\theta := \sM_\theta^\an$ carries a canonical
K\"ahler metric and $\cR_\rho := \sR_\rho^\an$ carries a natural
Hermitian connection.
In the two dimensional situation considered by Kronheimer and
Nakajima, the hyperk\"ahler condition ensures automatically that the
metric on $M_\theta$ is Ricci-f\/lat and the connection on $\cR_\rho$ is
HYM. In dimension three this is no longer true.
The main result of this paper shows that they can be deformed to satisfy these conditions and, most importantly, the resulting connection on $\cR = \bigoplus_\rho \cR_\rho$ is inf\/initesimally rigid.

\begin{Theorem}
  \label{Thm_KRF_HYM}
  Let $G$ be a finite subgroup of $\SL(3, \C)$ acting freely on $\C^3{\setminus}\{0\}$ and let $\theta\in \Theta_{\Q}$ be a generic rational stability parameter.
  Then the following hold:
  \begin{enumerate}\itemsep=0pt
  \item[$(1)$]
    $M_\theta$ carries a Ricci-flat ALE K\"ahler metric $g_{\theta, \rm RF}$.
  \item[$(2)$]
    For each $\rho \in \Irr(G)$ the tautological holomorphic line bundle $\cR_{\rho}$ carries an asymptotically flat HYM connection compatible with the holomorphic structure.
  \item[$(3)$]
    The induced HYM connection on the tautological bundle $\cR=\bigoplus_{\rho\in\Irr(G)}\cR_\rho$ is infinitesimally rigid.
  \end{enumerate}
\end{Theorem}

\begin{Remark}
  If $G$ does not act freely on $\C^3{\setminus} \{0\}$, then there are still crepant resolutions $\pi\colon \sM_\theta \to \C^3\!/G$;
  however, the asymptotic geometry of $M_\theta := \sM_\theta^\an$ will be quasi-asymptotically locally Euclidean (QALE).
  This causes a large number of additional dif\/f\/iculties. Some of these
  issues are tackled in recent work of the f\/irst named author~\cite{Degeratu2014}.
\end{Remark}

\begin{Remark}
  Theorem~\ref{Thm_KRF_HYM} is of interest in the context of higher dimensional gauge theory, and can be used,  for example, to extend the second named author's construction of $\rG_2$-instantons on generalised Kummer constructions~\cite{Walpuski2013} to $\rG_2$-manifolds arising from $\rG_2$-orbifolds with codimension $6$ singularities.
\end{Remark}

The existence of the Ricci-f\/lat K\"ahler metric on $M_{\theta}$ is a consequence of Joyce's proof of the Calabi conjecture for ALE crepant resolutions \cite[Section~8]{Joyce2000}, while the existence of the HYM connection is a simple consequence of the properties of the Laplace operator on ALE manifolds.
The most interesting and dif\/f\/icult part of Theorem~\ref{Thm_KRF_HYM} is the inf\/initesimal rigidity statement.
This will be a consequence of a vanishing result in
Lemma~\ref{lem:vanishing}, whose proof relies heavily on properties of the bounded derived category of coherent sheaves on $\sM_\theta$ and consequences of the derived category McKay correspondence.
This is in stark contrast with the work of Kronheimer and Nakajima~\cite{Kronheimer1990} who prove inf\/initesimal rigidity by bare hands.

By a result of Craw and Ishii, see Theorem~\ref{thm:ci21}, the tautological line bundles $\{\cR_{\rho}\colon \rho \in \Irr(G) \}$ form a basis in the $K$-theory of $M_\theta$ and thus their Chern characters form a basis of $H^*(M_{\theta}, \R)$. Since in our set-up $G$ acts freely on $\C^3{\setminus} \{0\}$, the exceptional divisor is contained in a compact subset of $M_\theta$ whose complement is homeomorphic to a truncated cone over $S^5/G$ and therefore $H^2(M_{\theta}, \R)\iso H^2_c(M_{\theta},\R)$. Hence we have the triple product
\begin{gather}\label{eq:triple-prod}
  \int_{M_\theta}\colon H^2(M_\theta,\R)^{\otimes 3} \iso H^2_c(M_\theta,\R)^{\otimes 3} \to\R.
\end{gather}
Exploiting the inf\/initesimal rigidity and using the Atiyah--Patodi--Singer index theorem applied to carefully chosen Dirac
operators we derive the following result for this triple product.

\begin{Theorem}
  \label{Thm_MultiplicativeIdentity}
  Let $G$ be a finite subgroup of $\SL(3, \C)$ acting freely on  $\C^3{\setminus}\{0\}$.
  Then for every generic rational stability parameter $\theta\in \Theta_\Q$, we have
  \begin{gather}\label{Eq_MultiplicativeIdentity}
   \frac12 \int_{M_\theta} c_1(\cR_{\rho})^2 c_1(\cR_{\sigma}) -
    c_1(\cR_{\rho}) c_1(\cR_{\sigma})^2 = - \big(C^{-1}\big)_{\rho\sigma}
  \end{gather}
  for all $\rho, \sigma \in \Irr_0(G)$.
  Here $C$ is a matrix which depends only on $G\subset \SL(3,\C)$, see~\eqref{Eq_C}.
\end{Theorem}

\begin{Remark}
  If we denote by $\wtilde\ch
  := \ch-\rk$  the reduced Chern character, then~\eqref{Eq_MultiplicativeIdentity} can equivalently be written as
  \begin{gather}\label{Eq_MultiplicativeIdentity2}
    \int_{M_\theta} \widetilde\ch(\cR_\rho)\widetilde\ch(\cR_\sigma^*) = -\big(C^{-1}\big)_{\rho\sigma},
  \end{gather}
  which is the natural analogue of Kronheimer and
  Nakajima's~\cite[Theorem~A.7]{Kronheimer1990} giving their
  geometrical interpretation of the McKay correspondence.
\end{Remark}

It has been pointed out to us by a referee of an earlier version of
this paper, that the formula~\eqref{Eq_MultiplicativeIdentity} could
also be derived from the work of Ito and Nakajima~\cite[Corollary~5.3]{Ito2000} by
applying the Riemann--Roch theorem for quasi-projective
varieties~\cite{Baum1979}.

Note that Theorem~\ref{Thm_MultiplicativeIdentity} exhibits a certain part of this triple product~\eqref{eq:triple-prod} that does not depend on the specif\/ic choice of crepant resolution but only on the subgroup~$G$ of~$\SL(3,\C)$.
It is now interesting to ask precisely how much of~\eqref{eq:triple-prod} is independent of the choice of crepant resolution and to try to determine~\eqref{eq:triple-prod} depending on $\theta$.
To our knowledge both of these question are still open and deserve to be investigated in future work.

The paper is organised as follows:
In Section~\ref{Sec_ModuliSpacesOfGConstellations} we brief\/ly recall the construction of crepant resolutions as moduli spaces of $G$-constellations, introduce the Fourier--Mukai transform and collect the results of Bridgeland, King and Reid~\cite{Bridgeland2001} and Craw and Ishii~\cite{Craw2004} that are
relevant for our work.
In Section~\ref{sec:kq} we present the construction of moduli spaces of $G$-constellations from the K\"ahler point of view and discuss their geometry in more detail.
In Section~\ref{sec:RFHYM} we prove the f\/irst two parts of Theorem~\ref{Thm_KRF_HYM}, while in Section~\ref{sec:rigid} we prove the inf\/initesimal rigidity statement.
Section~\ref{sec:dirac} introduces the Dirac operators on $M_\theta$ relevant for the proof of Theorem~\ref{Thm_MultiplicativeIdentity}, establishes their main properties, and uses the build-up of work so far to prove the vanishing of the index in Proposition~\ref{prop:mv}.
We complete the proof of Theorem~\ref{Thm_MultiplicativeIdentity} in Section~\ref{sec:mckay}.

\section[Moduli spaces of $G$-constellations]{Moduli spaces of $\boldsymbol{G}$-constellations}
\label{Sec_ModuliSpacesOfGConstellations}

Let $G$ be a f\/inite subgroup of $\SL(3,\C)$.  We denote by $\Irr(G)$
its set of irreducible representations, by~$\Rep(G)$ its
representation ring, and by~$R$ its regular representation.  Moreover,
$G$ has a natural action on~$\C^3$ which we tacitly assume
throughout this article.

\begin{Definition}
  \label{Def_GConstellation}
  A \emph{$G$-sheaf} on $\C^3$ is a coherent sheaf $\sF$ together
  with an action of $G$ which is equivariant with respect to the
  action of~$G$ on~$\C^3$.  A $G$-sheaf is called a
  \emph{$G$-constellation} if $H^{0}(\C^3, \sF) \iso R$ as
  $G$-modules.  Two $G$-constellations are \emph{isomorphic} if they
  are isomorphic as $G$-sheaves.
\end{Definition}

From this def\/inition it follows that the set-theoretic support of a $G$-constellation is a f\/inite union of $G$-orbits.
Thus a $G$-constellation is a sheaf-theoretic generalisation of the notion of $G$-orbit.

\begin{Definition}\label{Def_Theta}
  The set
  \begin{gather*}
    \Theta := \left\{ \theta\in\Hom_\Z(\Rep(G),\Z)\colon
                      \theta(R)=0 \right\}
  \end{gather*}
  is called \looseness= 1 the \emph{space of integral stability parameters}.  The
  sets $\Theta_\Q := \Theta \otimes_\Z \Q$ and $\Theta_\R := \Theta
  \otimes_\Z \R$ are~cal\-led the \emph{space of rational stability
    parameters} and the \emph{space of real stability parameters},
  respectively.  Given $\theta \in \Theta_\R$, a~$G$-constellation
  $\sF$ is called \emph{$\theta$-stable}
  (resp.\ \emph{$\theta$-semi-stable}) if each non-tri\-vial proper
  $G$-equivariant subsheaf $\sE\subset\sF$ satisf\/ies
  $\theta(H^0(\C^3, \sE))>0$ (resp.\ $\theta(H^0(\C^3, \sE))\geq
  0$).
\end{Definition}

When $\theta$ is a rational stability parameter GIT techniques are
used to prove that the $\theta$-stable $G$-constellations form a f\/ine
moduli space $\sM_{\theta}$.

\begin{Theorem}[Craw and Ishii \protect{\cite[Section~2.1]{Craw2004}}]\label{thm:ci21}
  If $\theta\in\Theta_\Q$, then there exists a fine moduli
  spa\-ce~$\sM_\theta$ of $\theta$-stable $G$-constellations on
  $\C^3$.  Moreover, for each representation~$\rho$ of~$G$ there
  exists a~locally free sheaf~$\sR_\rho$ on~$\sM_\theta$. If $\rho$
  and $\sigma$ are two representations of~$G$, then
  $\sR_{\rho\oplus\sigma}=\sR_\rho\oplus\sR_\sigma$.
\end{Theorem}

\begin{proof}[Sketch of the proof]
  The construction of $\sM_\theta$ is based on ideas of King
  \cite{King1994} and Sardo~In\-f\/irri~\cite{SardoInfirri1996}.  A
  $G$-constellation on $\C^3$ is a $G$-equivariant
  $\Sym^\bullet(\C^3)^*$-module structure on $R$, i.e., a~$G$-equivariant homomorphism $\Sym^\bullet(\C^3)^* \to \End(R)$.
  Hence, each point $B$ in
  \begin{gather}
    \label{eq:N}
    N := \big\{ B\in\big(\End(R)\otimes \C^3\big)^G \colon [B\wedge B]=0\in\End(R)\otimes\Lambda^2\C^3 \big\}
  \end{gather}
  def\/ines a $G$-constellation via $p \in
  \Sym^\bullet(\C^3)^* \mapsto p(B) \in \End(R)$.
  (Here $[\cdot \wedge \cdot]$ is composed of the commutator $[\cdot, \cdot]\colon \End(R)\otimes\End(R) \to \End(R)$ and the wedge product $\wedge\colon \C^3\otimes \C^3 \to \Lambda^2 \C^3$.)
  In fact, every
  $G$-constellation arises in this way.  Two points in~$N$ yield
  isomorphic $G$-constellations if and only if they are related by a
  $G$-equivariant automorphism of~$R$, i.e., an element of
  $\GL(R)^G$.  Since $R=\bigoplus_{\rho \in \Irr(G)} \C^{\dim\rho}
  \otimes \rho$, Schur's lemma gives
  \begin{gather*}
    \GL(R)^G=\prod_{\rho\in \Irr(G)} \GL\big(\C^{\dim\rho}\big).
  \end{gather*}
  Because the diagonal $\C^*\subset \GL(R)^G$ acts trivially on~$N$, the
  action of $\GL(R)^G$ descends to an action of $\PGL(R)^G$.

  An integral stability parameter $\theta\in\Theta$ determines a~character $\chi_\theta\colon \PGL(R)^G\to\C^*$ def\/ined by
  \begin{gather}\label{eq:chitheta}
    \chi_\theta([g])
      = \chi_\theta([(g_\rho)])
     := \prod_{\rho\in \Irr(G)} \det(g_\rho)^{\theta(\rho)}.
  \end{gather}
  King \cite[Proposition~3.1]{King1994} proved that an element of~$N$
  is stable (resp.~semi-stable) in the sense of GIT with respect to
  $\chi_\theta$ if and only if the corresponding $G$-constellation is
  $\theta$-stable (resp.\ $\theta$-semi-stable).  Let~$N^s_\theta$
  (resp.~$N^{ss}_\theta$) be the subset of GIT (semi-)stable points with
  respect to $\chi_\theta$ in $N$ and let
  \begin{gather*}
    \sM_\theta:=N^s_\theta/\PGL(R)^G
  \end{gather*}
  be the corresponding GIT quotient.  As schemes,
  $\sM_{k\theta}=\sM_{\theta}$ for any $k \in \N$; therefore, the
  above construction extends to rational stability parameters
  $\theta\in\Theta_\Q$ as well\footnote{For
    each $\theta\in\Theta_\Q$, we can f\/ind $k \in \N$ so that $k\theta \in \Theta$.
    We set $\sM_{\theta} := \sM_{k\theta}$.
    By what was said earlier, this does only depends on $\theta$.}.

  To see that $\sM_{\theta}$ is indeed a f\/ine moduli space of
  $\theta$-stable $G$-constellations, we construct a~\emph{universal
    $G$-constellation} $\sU_\theta$ on $\sM_\theta\times \C^3$.  For
  this purpose we identify
  \begin{gather}\label{eq:pglid}
    \PGL(R)^G \iso \prod_{\rho \in \Irr_0(G)}\GL\big(\C^{\dim\rho}\big),
  \end{gather}
  where $\Irr_0(G)$ is the set of non-trivial irreducible
  representations of $G$.  In this way $\PGL(R)^G$ acts on $R$.  This
  makes $R\otimes\sO_N$ into a $\PGL(R)^G$-equivariant sheaf on $N$.
  We denote its descend to $\sM_\theta$ by $\sR$.  Since the universal
  morphism $R\otimes\sO_N \to \C^3\otimes R\otimes\sO_N$ is
  $\PGL(R)^G$-equivariant, it descends to a universal morphism $\sR
  \to \C^3\otimes \sR$ on $\sM_{\theta}$. This determines the
  universal $G$-constella\-tion~$\sU_{\theta}$ on~$\sM_\theta\times\C^3$.
  Concretely, $\sU_\theta$ is the sheaf obtained by pulling back $\sR$ to $\sM_\theta\times\C^3$ (via the projection to $\sM_\theta$) with the action of $\sO_{\C^3}=\Sym^\bullet(\C^3)^*$ prescribed by the universal morphism.

  Let $\rho\colon G \to \Aut(R_{\rho})$ be a representation of $G$.  Then
  $\PGL(R)^G$ acts on $R_\rho$ via the identif\/ication~\eqref{eq:pglid}
  and as above we can associate with $\rho$ a locally free sheaf
  $\sR_{\rho}$ on $\sM_\theta$.  It is clear from the construction
  that $\sR_{\rho\oplus\sigma}=\sR_\rho\oplus\sR_\sigma$.
\end{proof}

\begin{Remark}\label{rem:map}
  When a $G$-constellation is $\theta$-stable for some real
  stability parameter $\theta \in \Theta_\R$, its set-theoretic
  support is a unique $G$-orbit of $\C^3$. Thinking of this as a
  point in $\C^3/G$, we obtain a well-def\/ined map $\pi_{\theta}\colon
  \sM_{\theta} \to \C^3/G$.
\end{Remark}

To obtain further insight into the spaces $\sM_{\theta}$ and the
properties of the map $\pi_{\theta}$, it is helpful to use the
language of derived categories.  We f\/irst recall the \emph{bounded
  derived category} $D(\sA)$ associated with an abelian category
$\sA$.  For details we refer the reader to B\"uhler's notes
\cite{Buehler2007} as well as Thomas' article \cite{Thomas2001} and
Huybrechts' book \cite{Huybrechts2006}, both of which underline the
importance of derived categories of coherent sheaves in algebraic
geometry.  Roughly speaking, $D(\sA)$ is obtained from the category of
bounded chain complexes in $\sA$ by formally inverting
quasi-isomorphisms.  If $A,B\in\sA$ are considered as bounded chain
complexes concentrated in degree zero, then $\Hom_{D(\sA)}(A,B)$ is a
complex whose cohomology computes $\Ext^\bullet(A,B)$, that is,
\begin{gather*}
  H^\bullet(\Hom_{D(\sA)}(A,B))=\Ext^\bullet(A,B).
\end{gather*}
If $\sB$ \looseness=-1 is another abelian category and $f\colon \sA \to \sB$ is a
left-exact functor, then one assigns to it a~\emph{right derived
  functor} ${\mathbf R} f \colon D(\sA) \to D(\sB)$.  If $A\in\sA$ is considered
as a bounded chain complex concentrated in degree zero, then ${\mathbf R} f(A)$ is a complex which computes $R^\bullet f(A)$, that is,
\begin{gather*}
  H^\bullet({\mathbf R} f(A))=R^\bullet f(A).
\end{gather*}
An analogous construction assigns to every right-exact functor $g
\colon \sA \to \sB$ a \emph{left derived functor}~$\bL g$.  As is
customary when working with derived categories, we will often write~$f$ and $g$ instead of~${\mathbf R} f$ and~$\bL g$.

An important example of a derived category is $D(\Coh(X))$, the
bounded derived category of coherent sheaves over a~scheme~$X$.  If
$X$ and $Y$ are two schemes and $K\in\Coh(X\times Y)$ is a~co\-he\-rent
sheaf, then the \emph{Fourier--Mukai transform with kernel~$K$} is the
functor $\Phi_K\colon D(\Coh(X))\to D(\Coh(Y))$ def\/ined by
\begin{gather*}
  \Phi_K(-):=(p_2)_*(p_1^*-\otimes K).
\end{gather*}
Here $p_1^*$, $(p_2)_*$ and $\otimes$ are taken in the derived sense,
with $p_1$ and $p_2$ denoting the projections from $X\times Y$ to $X$
and $Y$ respectively. A simple instance of a Fourier--Mukai transform
is the following: If $f\colon X\to Y$ is a morphism and $\sO_\Gamma$
denotes the structure sheaf of its graph $\Gamma\subset X\times Y$,
then~$\Phi_{\sO_\Gamma}$ is nothing but~$f_*$.

In our context, we denote by $D(\sM_\theta)$ the bounded derived
category of coherent sheaves on~$\sM_\theta$ and by $D^G(\C^3)$ the
bounded derived category of $G$-sheaves on $\C^3$, which is the same
as the bounded derived category $D([\C^3/G])$ of coherent sheaves on
the stack $[\C^3/G]$.  One of the key ideas of Bridgeland, King and
Reid~\cite{Bridgeland2001} was to introduce the Fourier--Mukai
transform $\Phi_\theta\colon D(\sM_\theta) \to D^G(\C^3)$ whose kernel is
given by the universal $G$-constellation $\sU_{\theta}$
\begin{gather*}
  \Phi_\theta(-)=q_*(p^*(-\otimes\rho_0)\otimes\sU_\theta)
\end{gather*}
to the study of crepant resolutions.  Here $p\colon \sM_\theta\times
\C^3\to\sM_\theta$ and $q\colon \sM_\theta\times \C^3\to \C^3$ are the
canonical projections and~$\rho_0$ is the trivial representation of~$G$.

\begin{Definition}
  \label{Def_ThetaGeneric}
  A real stability parameter $\theta\in\Theta_\R$ is called
  \emph{generic}, if there exists no non-trivial proper
  subrepresentation $S\subset R$ such that $\theta(S)=0$.
\end{Definition}

The space of generic stability parameters is dense in
$\Theta_\R$. Moreover, if $\theta \in \Theta_{\R}$ is generic, then
every $\theta$-semi-stable $G$-constellation is $\theta$-stable
(see also \cite[Section~2.2]{Craw2004}). In the particular case when~$\theta$ is a generic {\em rational} stability parameter, the above
techniques are used to prove that~$\sM_{\theta}$ together with the map
$\pi_{\theta}$ def\/ined in Remark~\ref{rem:map} is a crepant resolution
of singularities of~$\C^3/G$.

\begin{Theorem}[Craw and Ishii \protect{\cite[Proposition~2.2 and Theorem~2.5]{Craw2004}}]\label{thm:bkr}
  For each $\theta \in \Theta_{\Q}$ generic, the morphism
  $\pi_{\theta}\colon \sM_\theta\to\C^3/G$ is a projective crepant
  resolution and the Fourier--Mukai transform~$\Phi_\theta$ is an
  equivalence of derived categories.  Moreover, the locally free
  sheaves~$\sR_{\rho}$ form a~$\Z$-basis in $K$-theory.
\end{Theorem}

\begin{Remark}
  Bridgeland, King and Reid~\cite{Bridgeland2001} f\/irst proved this
  result for Nakamura's $G$-Hilbert scheme.  Craw and Ishii observed
  that their proof works more generally for moduli spaces of
  $G$-constellations.  In the course of the proof of the fact that~$\Phi_\theta$ is an equivalence of derived categories one has to
  show that $\sM_\theta$ is smooth.  This is achieved by appealing to
  a deep result from commutative algebra called the intersection
  theorem \cite[Theorem~7.1]{Bridgeland2001}.  That $\pi_\theta$ is a~crepant resolution then follows from a categorical criterion for a~resolution to be crepant \cite[Lemma~3.1]{Bridgeland2001}.
\end{Remark}

Furthermore, for abelian subgroups $G$ of $\SL(3,\C)$ toric geometry
techniques are used to prove a partial converse of
Theorem~\ref{thm:bkr}.

\begin{Theorem}[Craw and Ishii \protect{\cite[Theorem~1.1]{Craw2004}}]\label{thm:ci-2}
  If $G$ is an abelian subgroup of $\SL(3,\C)$, then every projective
  crepant resolution of $\C^3/G$ is a moduli space of $\theta$-stable
  $G$-constellations for some generic $\theta\in\Theta_{\Q}$.
\end{Theorem}

\section[$\sM_\theta$ via K\"ahler reduction]{$\boldsymbol{\sM_\theta}$ via K\"ahler reduction}
\label{sec:kq}

We now approach the previous discussion from the K\"ahler point of view.
There is no loss in assuming that the f\/inite group $G\subset\SL(3,\C)$
preserves the standard Hermitian metric on~$\C^3$, that is,
$G\subset\SU(3)$.  We also f\/ix a $G$-invariant Hermitian metric on
$R$.
Then the vector space $(\End(R)\otimes\C^3)^G$ naturally is a K\"ahler
manifold with K\"ahler form
\begin{gather*}
  \omega(B,C)
  := \operatorname{Im} \sum_{\alpha=1}^3 \tr(B_\alpha C_\alpha^*)
   = \sum_{\alpha=1}^3 \tfrac{1}{2i}\tr(B_\alpha C_\alpha^*-B_\alpha^*C_\alpha).
\end{gather*}
Here we identify $B\in(\End(R)\otimes\C^3)^G$ with a triple
$(B_1,B_2,B_3)$ of endomorphisms of $R$.

\begin{Proposition}
  The action of $\PU(R)^G$ on $(\End(R)\otimes\C^3)^G$ by conjugation
  is Hamiltonian with moment map $\mu\colon (\End(R)\otimes\C^3)^G\to
  (\fpu(R)^G)^*$ given by
  \begin{gather*}
    \langle \mu(B),\xi\rangle =\sum_\alpha \tfrac{1}{2i} \tr(\xi
    [B_\alpha,B_\alpha^*]).
  \end{gather*}
\end{Proposition}

\begin{proof}
  It is enough to prove this for the action of $U(R)^G$.  If
  $\xi\in\fu(R)^G$, then the corresponding vector f\/ield $X_\xi$ on
  $(\End(R)\otimes\C^3)^G$ is given by $X_\xi(B) = [\xi,B].$ Thus
  \begin{gather*}
    i(X_\xi)\omega(\hat B)
      =\sum_\alpha \tfrac{1}{2i}\tr\big([\xi,B_\alpha] \hat
     B_\alpha^*-[\xi,B_\alpha]^*\hat B_\alpha\big) \\
\hphantom{i(X_\xi)\omega(\hat B)}{}
=\sum_\alpha \tfrac{1}{2i}\tr\big(\xi\big([B_\alpha,\hat
     B_\alpha^*]+[\hat B_\alpha,B_\alpha^*]\big)\big)
     =\langle \rd\mu(B)\hat B,\xi\rangle .\tag*{\qed}
  \end{gather*}
  \renewcommand{\qed}{}
\end{proof}

To continue, we f\/irst need to analyze the relation between
$(\fpu(R)^G)^*$ and the spaces of stability parameters introduced in
Def\/inition~\ref{Def_Theta}.  For each $\theta\in\Theta_\R$, we def\/ine
$\zeta_\theta\in(\fpu(R)^G)^*$ by
\begin{gather*}
  \zeta_\theta(\xi)
    := -\sum_{\rho\in\Irr(G)} i\theta(\rho)\tr(\xi\cdot \pi_\rho)
\end{gather*}
for all $\xi \in \fpu(R)^G$. Here $\pi_{\rho}\colon R \to \C^{\dim \rho}
\otimes R_{\rho} $ is the projection onto the $\rho$-isotypical
component of the regular representation and $\xi \cdot \pi_{\rho}$
is thought of as an element in $\End(R)$. Since $\{i\pi_\rho\,|\, \rho
\in \Irr(G)\}$ spans the center of $\fu(R)^G$, we can identify
$\Theta_{\R}$ with the centre of $(\fpu(R)^G)^*$ via
$\theta\mapsto\zeta_\theta$. Under this identif\/ication, the generic
stability parameters $\theta \in \Theta_\R$ correspond to
$\zeta_{\theta}$ satisfying $\zeta_{\theta}(i\pi_S) \neq 0$ for all
non-trivial proper subrepresentations $S \subset R$. Here $\pi_S\colon R
\to S$ denotes the orthogonal projection onto $S$.  Moreover, when
$\theta$ is integral, then $\zeta_\theta=-i\rd \chi_{\theta} \in
(\fpu(R)^G)^*$ with~$\chi_{\theta}$ the character def\/ined in~\eqref{eq:chitheta}.

With the above identif\/ication, for each $\theta \in \Theta_\R$ we let
\begin{gather}\label{eq:Mtheta-K}
  M_\theta: = \big(N\cap \mu^{-1}(\zeta_{\theta})\big)\!/ \PU(R)^G
\end{gather}
be the corresponding K\"ahler quotient of the restriction of the moment
map $\mu$ to $N$ def\/ined in~\eqref{eq:N}. Note that $N$ is a complex
subvariety of $(\End(R) \otimes \C^3)^G$ and thus the K\"ahler quotient
makes sense.  The space $M_{\theta}$ comes with a natural set of
bundles constructed the following way: For each representation
$\rho\colon G \to \GL(R_{\rho})$ of $G$, the group~$\PU(R)^G$ acts on
$R_{\rho}$ via~\eqref{eq:pglid}.  Let
\begin{gather*}
  \cR_\rho:=\big(N\cap \mu^{-1}(\theta)\big)\times_{\PU(R)^G} R_{\rho}
\end{gather*}
be the associated complex vector bundle over $M_\theta$.  As we will
show in a moment, the bund\-les~$\cR_\rho$ carry natural holomorphic
structures.  We call $\cR_{\rho}$ the \emph{tautological $($holomorphic$)$
  bundle} associated with $\rho$.

We proceed now to describe the relation between $\sM_{\theta}$ and
$\sR_{\rho}$ def\/ined in Theorem~\ref{thm:ci21} and the $M_{\theta}$
and $\cR_{\rho}$ def\/ined above. Note that the f\/irst makes sense only
for the rational stability parameters $\theta$, while the second makes sense
for all $\theta\in\Theta_\R$.

\begin{Proposition}
  If $B\in N \cap \mu^{-1}(\zeta_\theta)$, then the $G$-constellation $\sF$
  associated to $B$ is $\theta$-semi-stable. Therefore we have
  $N\cap \mu^{-1}(\zeta_\theta) \subset N^{ss}_\theta$ for all $\theta \in
  \Theta_\R$.
\end{Proposition}

\begin{proof}
  Let $\sE$ be a non-trivial proper $G$-equivariant subsheaf of
  $\sF$.  Then the regular representation decomposes into two
  non-trivial proper subrepresentations $ R=S\oplus T$ with
  $S:=H^0(\C^3, \sE)$ and $T$ its orthogonal complement. Corresponding
  to $\sE$ there is an associated triple of matrices
  $C\in\End(S)\otimes \C^3$.  Moreover, since each component of $B$
  leaves $S$ invariant, $B= \left(\begin{smallmatrix} C & D \\ 0 &
      E \end{smallmatrix}\right)$ with $D\in\Hom(T,S)\otimes\C^3$ and $E \in \End(S)\otimes\C^3$.
  Therefore,
  \begin{gather*}
    \langle \mu(B),i\pi_S\rangle
    =\tfrac{1}{2}\tr_S([C,C^*]+DD^*)
    =\tfrac{1}{2}\tr_S(DD^*) \geq 0.
  \end{gather*}
  Since $\langle \mu(B), i \pi_S\rangle  = \zeta_\theta (i\pi_S) = \theta (S)$, it
  follows that $\theta(H^0(\C^3, \sE))\geq 0$.
\end{proof}

King \cite[Theorem~6.1]{King1994} proved the following version of the
Kempf--Ness theorem: If $\theta \in \Theta$, then each closed
$\PGL(R)^G$-orbit in $N^{ss}_\theta$ meets $N\cap
\mu^{-1}(\zeta_\theta)$ in precisely one $\PU(R)^G$-orbit.  From this
we obtain the following identif\/ication:

\begin{Proposition}\label{prop:git-kahler}
  Suppose that $\theta\in\Theta_\Q$ is generic.  Then the inclusion
  $N\cap \mu^{-1}(\zeta_\theta)\into N^s_\theta=N^{ss}_\theta$ induces
  a biholomorphic map from $M_\theta$ to the analytification of
  $\sM_\theta$.  This map identifies the complex vector bundle
  $\cR_\rho$ with the complex vector bundle underlying the
  analytification of the locally free sheaf $\sR_\rho$.
\end{Proposition}

The identif\/ication with $\sR_\rho^\an$ equips $\cR_\rho$ with a
holomorphic structure. Moreover, by Theorem~\ref{thm:ci21} we have
\begin{gather*}
  \cR_{\rho\oplus\sigma}=\cR_\rho\oplus\cR_\sigma.
\end{gather*}
as holomorphic vector bundles.

The K\"ahler quotient constrution above induces a metric $g_{\theta}$
and a K\"ahler form $\omega_{\theta}$ on $M_{\theta}$ for each $\theta
\in \Theta_\R$.  We also have a canonical connection $A_{\theta}$ on
the $\PU(R)^G$-bundle $\mu^{-1}(\zeta_\theta)\to M_\theta$ whose
horizontal space at $B\in\mu^{-1}(\zeta_\theta)$ is given by the
orthogonal complement in $T_B\mu^{-1}(\zeta_\theta)$ of the
tangent space to $\PU(R)^G$-orbit through $B$.  In the case when $G$
is a f\/inite subgroup of $\SU(3)$ which acts freely on $\C^3{\setminus}
\{0\}$ a more precise description of the geometry of $M_{\theta}$ and
$\cR_{\rho}$ can be given.  For this, we f\/irst need to recall a number
of def\/initions.

\begin{Definition}\label{def:ale}
  Let $G$ be a f\/inite subgroup of $\SO(n)$ acting freely on
  $\R^n{\setminus} \{0\}$.  A Riemannian manifold $(X,g)$ is called an
  \emph{asymptotically locally Euclidean (ALE) manifold asymptotic to
    $\R^n/G$ to order $\tau >0$} if there exists a compact subset $K
  \subset X$ and a dif\/feomorphism $\pi\colon (\R^n {\setminus} \bar
  B_1)/G \to X{\setminus} K$ so that
  \begin{gather*}
    |\del_k (\pi^*g-g_0)|_{g_0}=O\big(r^{-\tau-k}\big)
  \end{gather*}
  for all $0 \leq k \leq 2$.  Here we use the notation $r:=|x|$. In
  the above situation we also say that the metric $g$ is
  \emph{asymptotically locally Euclidean $($ALE$)$ of order~$\tau$}.
\end{Definition}

\begin{Definition}\label{def:afc}
  Let $H$ be a Lie group.  A connection $A$ on a $H$-bundle $E$ over
  an ALE manifold $(X,g)$ asymptotic to $\R^n/G$ is called
  \emph{asymptotically flat of order $\tau >0$} if there exists a f\/lat
  connection $A_0$ on $E|_{X{\setminus} K}$ such that
  \begin{gather*}
    \big|\nabla_{A_0}^k (A-A_0)\big| = O\big(r^{-\tau-k}\big)
  \end{gather*}
  for $0 \leq k \leq 1$.  Here we use a metric which is induced by the
  Euclidean metric on $\R^n$ and a~met\-ric on the adjoint bundle
  associated with $E$ which is parallel with respect to $A_0$.
\end{Definition}

With these def\/initions we can now characterise the geometry of
$M_{\theta}$ and of the corresponding tautological bundles
$\cR_{\rho}$.

\begin{Theorem}\label{thm:asympt}
  Let $G$ be a finite subgroup of $\SU(3)$ acting freely on
  $\C^3{\setminus} \{0\}$.  Then the following hold:
  \begin{enumerate}\itemsep=0pt
  \item[$(1)$] $(M_0,g_0)$ is isometric to the orbifold $\C^3/G$ with the
    induced orbifold K\"ahler metric. The corresponding connection~$A_0$
    is flat.
  \item[$(2)$] 
   If $\theta \in \Theta_\Q$ is generic,
    then $M_{\theta}$ is smooth and the induced K\"ahler metric
    $g_{\theta}$ is ALE of order~$4$.
  \item[$(3)$] 
  If $\theta\in \Theta_\Q$ generic, then
    the $\PU(R)^G$-connection $A_{\theta}$ is asymptotically flat of
    order~$1$.  Its curvature decays like~$r^{-4}$ and is of type
    $(1,1)$.  Moreover, the induced connection $A_{\theta,\rho}$ on~$\cR_\rho$ is compatible with the unique holomorphic structure on~$\cR_\rho$ for each $\rho \in \Irr(G)$.
  \end{enumerate}
\end{Theorem}

In the case of f\/inite subgroups of $\SU(2)$ the analogous theorem was proven by Kronheimer \cite{Kronheimer1989}, Kronheimer and Nakajima~\cite{Kronheimer1990}, and Gocho and Nakajima \cite{Gocho1992}.
For the above theorem, the smoothness of the K\"ahler quotient $M_{\theta}$ for generic $\theta\in \Theta_{\Q}$ follows from the identif\/ication with the algebraic quotient $\sM_{\theta}$ provided by Proposition~\ref{prop:git-kahler} and the result of Theorem~\ref{thm:bkr}.
The f\/irst statement and the remaining part of the second were proved by Sardo~Inf\/irri \cite{SardoInfirri1996} by generalising Kronheimer's proof.
The proof of the f\/irst two parts of the third statement is a direct generalisation of the proof in \cite[Proposition~2.2]{Kronheimer1990} and of Gocho and Nakajima's argument.
For the third part, note that the condition that $G$ is a f\/inite subgroup of $\SU(3)$ which acts freely on $\C^3{\setminus} \{0\}$ implies that $G$ is an abelian subgroup, and in fact it must be cyclic of prime order.
As a consequence, its irreducible representations are one dimensional, and the corresponding bundles $\cR_{\rho}$ are holomorphic line bundles.
Since $H^1(M_\theta,\cO_{M_\theta})=H^{0,1}(M_\theta)=0$, the holomorphic structure on $\cR_\rho$ is unique for all $\rho \in \Irr(G)$.
Therefore it must be compatible with the connection induced by $A_{\theta}$ on $\cR_{\rho}$.

\begin{Remark}
  The above characterization uses the identif\/ication in
  Proposition~\ref{prop:git-kahler} and therefore the smoothness of
  $M_{\theta}$ can only be inferred for $\theta$ a generic {\em
    rational} stability parameter.  It seems reasonable to expect that
  $M_\theta$ is smooth for all generic \emph{real} stability
  parameters $\theta \in \Theta_\R$.  Because of the homogeneity of
  the moment map, this is certainly true for all $t\theta$ with
  $\theta \in \Theta_{\Q}$ generic and $t\in(0,\infty)$.  One can show
  that for generic $\theta\in\Theta_\R$ the action of $\PU(R)^G$ on
  $N \cap \mu^{-1}(\zeta_\theta)$ is free.  To conclude that
  $M_{\theta}$ is smooth, however, one still needs to show that
  $N\cap \mu^{-1}(\zeta_\theta)$ is contained in the smooth locus of
  $N$.
\end{Remark}

\begin{Remark}\label{rem:fg}
  As mentioned in the introduction, the classif\/ication of f\/inite subgroups of $\SL(3,\C)$ was initiated by  Blichfeldt~\cite{Blichfeldt1917}. He, however, missed two groups, which were found much later by Yau and Yu~\cite{Yau-Stephan1993}.  According to the complete classif\/ication, there are ten families of f\/inite subgroups of~$\SL(3,\C)$.  The f\/irst family consists of abelian groups acting diagonally on~$\C^3$; all the other families contain only non-abelian groups. Direct examination shows that these
non-abelian groups do not act freely on~$\C^3{\setminus} \{0\}$. Using the structure of f\/inite abelian groups, it can be easily seen, that
$G$ can only act freely on $\C^3{\setminus} \{0\}$ if it is cyclic. Examining this, it follows that if a f\/inite subgroup of
$\SL(3,\C)$ acts freely on~$\C^3{\setminus} \{0\}$ then it must be isomorphic to~$\Z_n$ with~$n$ odd.
This statement, however, depends on the embedding on~$\Z_n$ into $\SL(3, \C)$. For example, $\Z_9$ can be embedded into $\SL(3,\C)$ in at least two ways: as the subgroup generated by the diagonal matrix $\operatorname{diag}(\epsilon,\epsilon^4, \epsilon^4)$, or as the subgroup generated by $\operatorname{diag}(\epsilon, \epsilon^3, \epsilon^5)$, with~$\epsilon$ a~$9^{\text{th}}$ root of unity. In the f\/irst instance $\Z_9$ acts freely on $\C^3{\setminus} \{0\}$, in the second, it does not.
\end{Remark}

\section[Ricci-f\/lat metrics on $M_\theta$ and HYM connections on $\cR_{\rho}$]{Ricci-f\/lat metrics on $\boldsymbol{M_\theta}$ and HYM connections on $\boldsymbol{\cR_{\rho}}$}
\label{sec:RFHYM}

The results of Kronheimer~\cite{Kronheimer1989} and of Gocho and Nakajima~\cite{Gocho1992} for f\/inite subgroups of $\SU(2)$ we alluded to above actually establish that the metric $g_{\theta}$ is Ricci-f\/lat and the induced connections on $\cR_{\rho}$ are anti-self-dual and, hence, Hermitian--Yang--Mills (HYM).
This is because in dimension two~$M_\theta$ arises via hyperk\"ahler reduction with the subspace $N$ being the zero locus of the complex component of the hyperk\"ahler moment map.
In higher dimension this is no longer the case.
Consequently, the metric $g_\theta$ on $M_\theta$ given by
Theorem~\ref{thm:asympt} is not necessarily Ricci-f\/lat and the
connections $A_{\theta,\rho}$ on~$\cR_\rho$ are not necessarily HYM.
Indeed, Sardo~Inf\/irri \cite[Example~7.1]{SardoInfirri1996} showed
that for $G=\Z_3\iso\langle \diag(e^{2\pi i/3},e^{2\pi i/3},e^{2\pi
  i/3})\rangle \subset\SU(3)$ and an appropriate choice of a generic
stability parameter, the corresponding K\"ahler quotient is the total
space of the line bundle $\cO_{\P^2}(-3)$ with a metric which has
non-vanishing Ricci curvature.
In this section we show that $M_\theta$ does admit a Ricci-f\/lat
K\"ahler metric and that the tautological line bundles $\cR_\rho$ carry asymptotically f\/lat HYM
connections, thus proving the f\/irst two parts of Theorem~\ref{Thm_KRF_HYM}.

The existence of the Ricci-f\/lat K\"ahler metric follows from the
following result:

\begin{Theorem}[Joyce \protect{\cite[Theorem~8.2.3]{Joyce2000}}]\label{thm:rf}
  Let $G$ be a finite subgroup of $\SU(n)$ acting freely on
  $\C^n{\setminus} \{0\}$.  Let $X$ be a smooth crepant resolution of
  $\C^n/G$ with an ALE K\"ahler metric $g$ of order $\tau>n$.  Then
  there exists a unique Ricci-flat ALE K\"ahler metric $g_\RF$ in the
  K\"ahler class of $g$.  The metric $g_\RF$ is ALE of order $2n$.
\end{Theorem}

\begin{Remark}
  Joyce states this result only for ALE K\"ahler metrics of order $\tau
  =2n$; however, his proof goes through for $\tau >n$. More specif\/ically, we need to modify the metric $g$ in its K\"ahler class to be f\/lat outside a compact set and then apply Joyce's proof of the Calabi conjecture in ALE set-up. For this modif\/ication of the metric, we need the $dd^c$-lemma to hold for a certain exact, real $(1,1)$-form which decays like $\rho^{-\tau}$ on the ALE end.
In~\cite[Theorem~8.4.4]{Joyce2000} this is done via a Stokes' theorem argument. To use this argument and conclude that the boundary integral is zero, one needs that $\tau > n$.
\end{Remark}

By Theorem~\ref{thm:asympt}(2) the induced K\"ahler
metric on $M_\theta$ is ALE of order $4$ and, hence, the above theorem
applies.

\begin{Corollary}
  Let $\theta \in \Theta_{\Q}$ be generic and let
  $(M_{\theta},g_{\theta})$ be the corresponding K\"ahler quotient. Then
  there exists a Ricci-flat ALE K\"ahler metric $g_{\theta, \RF}$ of
  order $6$ on $M_\theta$ in the same K\"ahler class as $g_\theta$.
\end{Corollary}

\begin{Remark}
  Note that in dimension $n\geq 4$, Theorem~\ref{thm:rf} does not
  apply anymore for the ALE K\"ahler metrics constructed on crepant
  resolutions of $\C^n/G$ via the K\"ahler reduction~\eqref{eq:Mtheta-K}
  since the decay $\tau=4$ is now too weak. Another argument is then needed to show
  the existence of Ricci-f\/lat ALE K\"ahler metrics on these crepant
  resolutions.
\end{Remark}

We now proceed to show the existence of asymptotically f\/lat HYM
connections on the tautological line bundle $\sR_\rho$.  Recall that a
$(1,1)$-connection on a complex vector bundle over a K\"ahler manifold
is called {\em Hermitian--Yang--Mills $($HYM$)$} if the contraction of its
curvature with the K\"ahler form vanishes identically, i.e., $\Lambda F_A
=0$. It turns out that it is a little easier to prove the existence
result in terms of {\em Hermitian--Yang--Mills $($HYM$)$ metrics}. These
are Hermitian metrics on holomorphic bundles with the property that
their Chern connection, the unique metric connection associated to
the holomorphic bundle, is HYM.

\begin{Definition}\label{def:afm}
  Let $E$ be a complex vector bundle over an ALE manifold $(X,g)$, let
  $h_0$ be a Hermitian metric on $E|_{X{\setminus} K}$ and let $A_0$ be
  a connection on $E|_{X{\setminus} K}$ compatible with $h_0$.  A~Hermitian metric $h$ on $E$ is called \emph{asymptotic to~$h_0$ to
    order $\tau >0$} if
  \begin{gather*}
    |\nabla_{A_0}^k (h-h_0)|_{h_0} = O\big(r^{-\tau-k}\big)
  \end{gather*}
  for all $0 \leq k \leq 2$.
\end{Definition}

\begin{Proposition}\label{prop:hym}
  Let $X$ be an ALE K\"ahler manifold and let $\cL$ be a holomorphic
  line bundle over $X$.  If $h_0$ is a Hermitian metric on $\cL$ such
  that the curvature $F_{h_0}$ of the Chern connection on $\cL$
  compatible with $h_0$ satisfies
  \begin{gather*}
    \Lambda F_{h_0} = O\big(r^{-2-\eps}\big)
  \end{gather*}
  for some $\eps>0$, then there exists a HYM metric $h$ on $\cL$.
  Moreover, for every $\tau \in (0, \epsilon)$ this metric is the
  unique HYM metric asymptotic to $h_0$ to order $\tau$.
\end{Proposition}

From the construction of the bundles $\cR_{\rho}$ and of the
corresponding connections $A_{\theta,\rho}$ in
Theorem~\ref{thm:asympt}(3), we see that
$A_{\theta, \rho}$ is the Chern connection of a Hermitian metric on
$\cR_{\rho}$. Therefore $\cR_{\rho}$ with this Hermitian metric
satisf\/ies the conditions of the above proposition yielding the desired
existence result:

\begin{Corollary}
  Let $\theta \in \Theta_{\Q}$ be generic.  Then for each
  $\rho\in\Irr(G)$ the tautological line bundle~$\cR_\rho$ on
  $M_\theta$ carries a~HYM $\U(1)$-connection with respect to
  $g_{\theta,\RF}$ which is asymptotically flat of order $\tau \in (0,2)$.
\end{Corollary}

\begin{Remark}
  Using some of the results derived in Section~\ref{sec:rigid}, one
  can show that the HYM connection associated with $h$ in
  Proposition~\ref{prop:hym} is asymptotically f\/lat of order $5$ and,
  hence, the Hermitian metric $h$ is asymptotic to a f\/lat metric to
  order $4$.
\end{Remark}

\begin{Remark}
  Using heat f\/low methods, Bando \cite{Bando1993} proved that every
  holomorphic bundle $\cE$ over an ALE K\"ahler manifold which admits a
  Hermitian metric $h_0$ with $|F_{h_0}|=O(r^{-2-\eps})$ does in
  fact carry a HYM metric.
\end{Remark}

The case of line bundles is much simpler than Bando's result and the
proof of Proposition~\ref{prop:hym} follows from the fact that the
Laplace operator is an isomorphism between certain weighted Sobolev
spaces.  Concretely, let $(X, g)$ be an ALE manifold asymptotic to
$\R^n/G$ as def\/ined in Def\/inition~\ref{def:ale} and let $r\colon
X\to[1,\infty)$ denote a smooth extension of the radius function from
$X{\setminus} K \iso (\R^n{\setminus} \bar B_1)/G$ to all of $X$.  For a~non-negative integer $k$ and a real number $\delta$ we denote by
$W^{k,2}_{\delta}(X)$ the weighted Sobolev space obtained as the
completion of $C^\infty_0(X)$ with respect the norm
\begin{gather}\label{eq:wsf}
  \|f\|_{W^{k,2}_\delta}: = \sum_{j=0}^k \big\|r^{-\delta-n/2+j}\nabla ^jf\big\|_{L^2}.
\end{gather}
Let $\Delta_{\delta}\colon W^{k+2,2}_\delta(X) \to W^{k,2}_{\delta-2}(X)$
denote the corresponding completion of the Laplacian $\Delta$.

\begin{Proposition}[Bartnik~\protect{\cite[Proposition~2.2]{Bartnik1986}}]\label{prop:Laplace}
  For $\delta\in(-n+2,0)$ the operator $\Delta_{\delta}$ is an isomorphism.
\end{Proposition}

\begin{proof}[Sketch of the proof]
  The weighted Laplacian $\Delta_\delta \colon W^{k+2,2}_\delta(X) \to W^{k,2}_{\delta-2}(X)$ is a Fredholm operator if and only if the weight parameter $\delta$ is not contained in its set of indicial roots at inf\/inity.
  This is a discrete set of real numbers which does not intersect the
  interval $(-n+2, 0)$, see Bartnik \cite[Sections~1 and~2]{Bartnik1986} for details.  Moreover, for $\delta<0$ the kernel of
  $\Delta_\delta$ is trivial by the maximum principle.  On the other
  hand, the cokernel of $\Delta_{\delta}$ is isomorphic to the kernel
  of its formal adjoint $\Delta_{n-2-\delta}$.  Therefore,
  $\Delta_\delta$ is an isomorphism for $\delta \in (-n+2,0)$.
\end{proof}

\begin{proof}[Proof of Proposition \ref{prop:hym}]
  Any Hermitian metric on $\cL$ is of the form $h = e^f h_0$, for some
  $f \in C^\infty(X)$ and $F_h=F_{h_0} + \delbar\del f \in
  \Omega^2(X,i\R)$. Therefore,
  \begin{gather*}
    i\Lambda F_h = i\Lambda F_{h_0} + \tfrac12 \Delta f.
  \end{gather*}
  Since $\Lambda F_{h_0}\in L^2_{-2-\tau}(X):=W^{0,2}_{-2-\tau}(X)$ for
  every $\tau\in(0,\epsilon)$, by Proposition \ref{prop:Laplace},
  there exists a unique $f\in W^{2,2}_{-\tau}(X)$ such that $\Delta f
  =-2i\Lambda F_{h_0}$. Moreover, a computation using the explicit
  form of the Green function on the end of the ALE manifold $X$ gives
  that $f=O(r^{-{\tau}})$ for all $0<\tau<\eps$, cf.~\cite[(proof of) Theorem~8.3.5]{Joyce2000}.
  \end{proof}

\section[Rigidity of HYM connections on the holomorphic tautological bundles]{Rigidity of HYM connections\\ on the holomorphic tautological bundles}\label{sec:rigid}

In this section we prove the inf\/initesimal rigidity statement in
Theorem~\ref{Thm_KRF_HYM}(3).  This will be an immediate
consequence of the following lemma, which is the core vanishing result
of this paper.

\begin{Lemma}\label{lem:vanishing}
  Let $\theta\in \Theta_\Q$ be generic and let $M_\theta$ be equipped
  with an ALE K\"ahler metric $g$. Let~$h$ be a Hermitian metric on the
  holomorphic bundle $\cR = \bigoplus_{\rho\in \Irr(G)} \cR_{\rho}$
  whose associated Chern connection $A$ is asymptotically flat of
  order $\tau >0$. Then the space
  \begin{gather*}
    \cH^1_A
    :=\Big\{ a\in\Omega^{0,1}(M_\theta,\cEnd(\cR))\colon
      \delbar_A a=\delbar_A^* a=0 ~\text{and}~
      \lim_{r\to\infty}\sup_{\del B_r}|a|=0 \Big\}
  \end{gather*}
  is trivial.
\end{Lemma}

Note that if the connection $A$ is HYM, then $\cH^1_A$ is its space of
inf\/initesimal deformations.  In particular, it follows that the HYM
connection on~$\cR$ induced by then HYM connections on the bundles~$\cR_\rho$ constructed in the second part of Theorem~\ref{Thm_KRF_HYM} is
inf\/initesimally rigid. This thus completes the proof of the third part
of Theorem~\ref{Thm_KRF_HYM}.

The strategy for proving Lemma~\ref{lem:vanishing} is as follows: We
f\/irst reduce to a vanishing result in complex geometry, see
Propositions~\ref{prop:d} and~\ref{prop:av}. Then since $\theta$ is a
generic stability parameter, Proposition~\ref{prop:git-kahler} gives
that $M_{\theta}$ is the analytif\/ication of the moduli space of
$\theta$-stable $G$-constellations $\sM_{\theta}$. Using GAGA, we
translate the vanishing into an algebraic geometry problem,
see~\eqref{eq:h10}, which we then solve using the results of
Bridgeland, King and Reid \cite{Bridgeland2001} and Craw and Ishii
\cite{Craw2004} for the moduli spaces of $G$-constellations and the
corresponding tautological free sheaves discussed in Section
\ref{Sec_ModuliSpacesOfGConstellations}.

It is a useful heuristic to think of bundles with decaying connections
as bundles on a compactif\/ication whose restrictions to the ``divisor
at inf\/inity'' satisfy certain ``vanishing conditions''.  With this in
mind, we compactify $M_\theta$ at inf\/inity by gluing~$M_\theta$ and
$(\P^3{\setminus}\{[0:0:0:1]\})/G$ along
$M_\theta{\setminus}\pi_{\theta}^{-1}(0)=(\C^3{\setminus}\{0\})/G$.  The
resulting space $\bar M_\theta$ is not a~complex manifold, but rather
a~complex orbifold.  One can think of $\bar M_\theta$ as obtained from
$M_\theta$ by adjoining the divisor $D=\P^2/G$ at inf\/inity.  $D$~is a~smooth orbifold divisor, i.e., it lifts to a smooth divisor in covers
of the uniformising charts.  The holomorphic bundle $\cR$ extends over
$D$ to a holomorphic bundle~$\bar \cR$ on~$\bar M_\theta$.  The
following result reduces the proof of Lemma~\ref{lem:vanishing} to a~problem in complex geometry.

\begin{Proposition}\label{prop:d}
  $\cH^1_A$ injects into $H^1(\bar M_\theta,\cEnd(\bar \cR)(-D))$.
\end{Proposition}

Recall that for a holomorphic vector bundle $\cE$, $\cE(-D)$ is the
sheaf of holomorphic sections of~$\cE$ vanishing to f\/irst order along~$D$.  The proof of Proposition~\ref{prop:d} requires two preparatory
results.

\begin{Proposition}\label{prop:r}
  Let $Z$ be a complex orbifold, $D$ be a smooth divisor in~$Z$ and
  $\cE$ be a holomorphic bundle on~$Z$. Denote by $i\colon D\into Z$ the
  inclusion of $D$ into $Z$.  Then the complex of sheaves
  $(\cA^\bullet,\delbar)$ defined by
  \begin{gather*}
    \cA^k (U):=\big\{ \alpha\in\Omega^{0,k}(U,\cE)\colon
                       i^*\alpha=0 \big\}
  \end{gather*}
  for $U\subset Z$ open is an acyclic resolution of $\cE(-D)$.
\end{Proposition}

\begin{proof}
  Since $i^*$ and $\delbar$ commute, $\cA^\bullet$ forms a
  complex.  Moreover, it is clear that $\cE(-D)$ is the kernel of
  $\cA^{0} \stackrel{\delbar}{\to} \cA^{1}$.

  The proof that $(\cA^\bullet,\delbar)$ is a resolution uses two
  ingredients: the Grothendieck--Dolbeault lemma and the fact that if~$U$ is a suf\/f\/iciently small open set, then holomorphic sections on
  $D\cap U$ extend to~$U$. First we show that these assertions also
  hold for orbifolds.  Let $U$ be a small open set which is covered by
  a uniformising chart $\tilde U/\Gamma$ where $\Gamma$ is a~f\/inite
  group. Lifting everything up to $\tilde U$, $\cE$ corresponds to a
  $\Gamma$-equivariant holomorphic bundle~$\tilde\cE$ and~$D$ to a
  $\Gamma$-equivariant smooth divisor $\tilde D$.  If
  $\alpha\in\Omega^{0,k}(U,\cE)$ satisf\/ies $\delbar\alpha=0$, then so
  does its lift $\tilde\alpha\in\Omega^{0,k}(\tilde
  U,\tilde\cE)^\Gamma$.  If $U$ (and thus~$\tilde U$) is suf\/f\/iciently
  small, then the usual Grothendieck--Dolbeault lemma yields
  $\tilde\beta\in\Omega^{0,k-1}(\tilde U,\tilde\cE)$ satisfying
  $\delbar\tilde\beta=\tilde\alpha$.  There is no loss in assuming
  that $\tilde\beta$ is $\Gamma$-invariant and thus pushes down to
  the desired primitive $\beta\in\Omega^{0,k-1}(U,\cE)$ of $\alpha$.
  We thus obtain the Grothendieck--Dolbeault lemma for orbifolds.
  Now, if~$s$ is a holomorphic section of~$\cE$ over $D\cap U$, we
  lift it to the uniformising chart $\tilde U$, where, provided~$U$ is
  suf\/f\/iciently small, we f\/ind a~$\Gamma$-equivariant extension. We
  then push this extension down to~$U$.  This proves the second
  assertion.

  Let now $U$ be a small open set of $Z$ and let
  $\alpha\in\Omega^{0,k}(U,\cE)$ with $\delbar\alpha=0$.  By the
  Grothendieck--Dolbeault lemma after possibly shrinking~$U$, we can
  f\/ind $\beta\in\Omega^{0,k-1}(U,\cE)$ satisfying $\delbar\beta=\alpha$.
  If $k\geq 2$, we apply the Grothendieck--Dolbeault lemma once more to
  obtain $\gamma\in \Omega^{0,k-2}(U\cap D, \cE)$ such that
  $\delbar\gamma=i^*\beta$.  We extend $\gamma$ smoothly to all of~$U$. Then $\beta-\delbar\gamma \in \cA^{k-1}(U)$ yields the desired
  primitive of~$\alpha$ on~$U$.  If $k=1$, we know that $\beta$
  restricts to a~holomorphic section~$\beta|_D$ of~$\cE|_{U\cap D}$,
  which can be extended to a holomorphic section $\delta$ on~$U$.
  Hence, $\beta-\delta\in\cA^0(U)$ is the desired primitive of~$\alpha$.

  Finally, $(\cA^\bullet,d)$ is an acyclic resolution of $\cE(-D)$,
  since the sheaves $\cA^\bullet$ are $C^\infty$-modules and
  therefore soft.
\end{proof}

\begin{Remark}
  In the def\/inition of $\cA^k$ it is not strictly necessary to require
  that $\alpha$ be smooth.  In fact, a simple application of elliptic
  regularity shows that it suf\/f\/ices that elements of~$\cA^k$ be in
  the H\"older space $C^{n-k,\gamma}$, where~$n$ denotes the complex
  dimension of~$Z$.
\end{Remark}

\begin{Proposition}\label{prop:decay}
  If $a\in\cH^1_A$, then
  \begin{gather}\label{eq:decay}
    \nabla_A^k a=O\big(r^{-5-k}\big) \qquad \text{for all $k\geq 0$.}
  \end{gather}
\end{Proposition}

\begin{proof}
  First observe that using simple scaling considerations and standard
  elliptic theory, \eqref{eq:decay} for $k>0$ follows from the case
  $k=0$.

  It is rather straightforward to obtain $a=O(r^{-4})$ using the
  maximum principle.  To obtain the stronger decay estimate it is
  customary to use a ref\/ined Kato inequality, see, e.g., Bando, Kasue
  and Nakajima~\cite{Bando1989}.  Recall that the classical Kato
  inequality is a consequence of the Cauchy--Schwarz inequality
  $|\langle \nabla_A a,a\rangle |\leq |\nabla_A a|\ |a|$.  In our case, since $a$ is not arbitrary but satisf\/ies $\delbar a = \delbar_A^* a=0$, the above inequality can be improved upon:
  There exists a constant $\gamma<1$ such that $|\rd|a||\leq \gamma|\nabla_A a|$ on
  the set $U:=\{x \in M_{\theta}\colon a(x) \neq 0\}$.  A detailed
  analysis shows that since we are working on a~$6$-dimensional real
  manifold, we can choose $\gamma$ to be $\sqrt{5/6}$, see, e.g.,~\cite{Calderbank2000}.

  We set $\gamma := \sqrt{5/6}$ and let $\sigma
  :=2-1/\gamma^2=4/5$. Using the ref\/ined Kato inequality for~$a$, we
  obtain
  \begin{gather*}
    (2/\sigma) \Delta |a|^\sigma
      = |a|^{\sigma-2}\big(\Delta|a|^2-2(\sigma-2)|\rd|a||^2\big) \\
\hphantom{(2/\sigma) \Delta |a|^\sigma}{} \leq |a|^{\sigma-2}\big(\Delta|a|^2+2|\nabla_A a|^2\big)
     = |a|^{\sigma-2}\langle a,\nabla_A^*\nabla_Aa\rangle
   \end{gather*}
   on $U$. The Weitzenb\"ock formula for $\nabla_A^*\nabla_Aa$ gives
   \begin{gather*}
    (2/\sigma) \Delta |a|^\sigma
      \leq |a|^{\sigma-2}
     \big(\langle \Delta_{\delbar_A} a,a\rangle
       +\langle \{{\rm Riem},a\},a\rangle +\langle \{F_A,a\},a\rangle \big),
   \end{gather*}
   with ${\rm Riem}$ the Riemannian curvature and $F_A$ the
   curvature of the connection $A$.  Because~$\Delta_{\delbar_A} a=0$
   and since by hypothesis of Lemma~\ref{lem:vanishing} the metric on
   $M_{\theta}$ is ALE and the connection~$A$ is asymptotically f\/lat,
   there exist constants $c,\beta>0$ so that on~$U$ we have
   \begin{gather*}
     (2/\sigma) \Delta |a|^\sigma \leq c (1+r)^{-2-\beta} |a|^{\sigma}.
   \end{gather*}
   Set $f: =|a|^\sigma$.  We show that $f=O(r^{-4})$, which is
   equivalent to the desired decay estimate for~$a$.  Note that on
   $U$ the above estimate gives
   \begin{gather}\label{eq:d-1}
     \Delta f \leq \frac{cf}{1+r^{2+\beta}}.
   \end{gather}
   Since $f$ decays at inf\/inity and is thus bounded, using the method of proof of
   \cite[Theorem~8.3.6(a)]{Joyce2000}, there exists $g\in
   C^{2,\alpha}$ with $g=O(r^{-\beta})$ such that
    \begin{gather*}
     \Delta g
     =\begin{cases}
       (\Delta f)^+ & \text{on $U$}, \\
       0 & \text{on}~M_{\theta}{\setminus} U .
     \end{cases}
  \end{gather*}
  Here $h^+$ denotes the positive part of the function $h$, i.e., $h^+
  := \max \{h, 0\}$.  Since $g$ is superharmonic and decays to zero at
  inf\/inity, the maximum principle implies that $g$ is non-negative on
  the boundary of~$U$; hence, by the maximum principle $f\leq g =
  O(r^{-\beta})$.  By~\eqref{eq:d-1} we then have $(\Delta f)^{+}=
  O(r^{-2-2\beta})$, which then yields that $f=O(r^{-2\beta})$.
  Iterating this argument $k$ times we obtain $f=O(r^{-k\beta})$ for
  all $k< (n-2)/\beta$ with $n = 6$ the real dimension of $M_\theta$.
  For the biggest~$k$ with this property, we have $2+(k+1)\beta>n$.
  Then, by \cite[Theorem~8.3.6(b)]{Joyce2000}, we can chose~$g$ above
  such that $g=O(r^{-4})$.  Therefore, $f=O(r^{-4})$ as desired.
\end{proof}

With these two preliminary results, we can now prove Proposition \ref{prop:d}.

\begin{proof}[Proof of Proposition \ref{prop:d}]
  Given $a \in \cH^1_A$, we extend it to a~$1$-form on $\bar
  M_\theta$ vanishing along~$D$.  From Proposition~\ref{prop:decay} it
  follows that $a$ vanishes to third order along $D$.  Hence, $a$ is
  in $C^{2,\alpha}(\bar M_\theta)$ and we can regard it as an element
  of $\cA^1(\bar M_\theta)$.  Since $\delbar a=0$, by
  Proposition~\ref{prop:r} it gives an element $[a]\in H^1(\bar
  M_\theta, \cEnd(\bar\cR)(-D))$.  This def\/ines a linear map $i\colon
  \cH^1_A\to H^1(\bar M_\theta,\cEnd(\bar\cR)(-D))$.

  We will now show that $i$ is injective.  Suppose that $i(a)=0$,
  i.e., there exists $b \in \cA^{0}(\bar M_\theta)$ so that $a =
  \delbar b$.  Since $b$ vanishes along $D$, its restriction to
  $M_\theta$ decays like $r^{-1}$. Using this together with
  $a=O(r^{-5})$, we can integrate by parts to obtain
  \begin{gather*}
    \|a\|_{L^2}^2= \int_{M_\theta} \langle a,\delbar b\rangle  \dvol_g
                = \int_{M_\theta}  \langle \delbar_A^*a,b\rangle  \dvol_g
                =0.
  \end{gather*}
  It follows that $a$ vanishes, and thus $i$ is injective.
\end{proof}

To prove Lemma~\ref{lem:vanishing} it now suf\/f\/ices to establish the
following result:

\begin{Proposition}\label{prop:av}
  $H^1(\bar M_\theta,\cEnd(\bar\cR)(-D))=0$.
\end{Proposition}

To prove this statement, we convert it into a problem in algebraic
geometry.  In the same way we compactif\/ied $M_\theta$, we can complete
the scheme $\sM_\theta$ at inf\/inity by attaching $D = \P^2/G$.  This
yields an algebraic stack $\bar\fM_\theta$.  Moreover, $\sR$ extends
to a locally free sheaf $\bar\sR$ on $\bar\fM_\theta$.  By GAGA
\cite[Th\'eor\`eme 5.10]{Toen1999}, Proposition~\ref{prop:av} is
equivalent to
\begin{gather}\label{eq:h10}
  H^1(\bar \fM_\theta,\sEnd(\sR)(-D))=0.
\end{gather}
To establish this we need the following consequence of
Theorem~\ref{thm:bkr}.

\begin{Proposition}\label{prop:ci}
  For generic $\theta\in \Theta_\Q$,
  \begin{gather*}
    H^k(\sM_\theta,\sR_\rho^*\otimes\sR_\sigma)=H^k\big(\C^3,\sO\otimes
    R_\rho^* \otimes R_\sigma\big)^G,
  \end{gather*}
  for all $\rho, \sigma \in \Irr(G)$.
    In particular, for $k>0$,
  \begin{gather*}
    H^k(\sM_\theta,\sR_\rho^*\otimes\sR_\sigma)=0.
  \end{gather*}
  If $G$ acts freely on $\C^3{\setminus}\{0\}$, we have a commutative
  diagram
  \begin{equation*}
    \begin{tikzpicture}[baseline=(current  bounding  box.center)]
      \matrix (m) [matrix of math nodes,row sep=3em,column
      sep=3em,minimum width=2em ] {
        H^k(\sM_\theta,\sEnd(\sR)) &
        H^k\big(\sM_\theta {\setminus}{\pi_{\theta}^{-1}(0)},\sEnd(\sR)\big)
        \\
        H^k\big(\C^3,\sO\otimes\End(R)\big)^G  &
        H^k\big((\C^3{\setminus}\{0\})/G,\sO\otimes\End(R)\big), \\
      };
      \path[-stealth]
      (m-1-1) edge node [above] {$i^*$} (m-1-2)
      (m-2-1) edge node [below] {$j^*$} (m-2-2)
      (m-1-1) edge node [left] {$\Phi_\theta$} (m-2-1)
      (m-1-2) edge node [right] {$(\pi_\theta)_*$} (m-2-2);
    \end{tikzpicture}
  \end{equation*}
  where $i\colon \sM_\theta{\setminus}\pi_{\theta}^{-1}(0)\to\sM_\theta$ and
  $j\colon \C^3{\setminus}\{0\}\to\C^3$ are the inclusion maps.
\end{Proposition}

\begin{proof}
  The f\/irst part is due to Craw and Ishii \cite[Lemma~5.4]{Craw2004}.
  Let us brief\/ly recall their proof.  We have
  \begin{gather*}
    H^k(\sM_\theta,\sR_\rho^*\otimes\sR_\sigma)
    =\Ext^k(\sO,\sR_\rho^*\otimes\sR_\sigma)
     =\Ext^k(\sR_\rho,\sR_\sigma)
    =H^k(\Hom_{D(\sM_{\theta})}(\sR_\rho,\sR_\sigma))
  \end{gather*}
  and
  \begin{gather*}
    H^k\big(\C^3,\sO\otimes R_\rho^* \otimes R_\sigma\big)^G
     =G - \Ext^k(\sO \otimes R_\rho,\sO\otimes R_\sigma)\\
     \hphantom{H^k\big(\C^3,\sO\otimes R_\rho^* \otimes R_\sigma\big)^G}{}
     =H^k\big(\Hom_{D^G(\C^3)}(\sO \otimes R_\rho,\sO\otimes R_\sigma)\big).
  \end{gather*}
  Moreover, the inverse of the Fourier--Mukai transform $\Phi_\theta$
  is given by
  \begin{gather*}
    \Phi_\theta^{-1}
    =\big(p_*\big(q^*(-)\otimes\sU_\theta^D[3]\big)\big)^G
    =\Big({-}\otimes \bigoplus_\rho \sR_\rho^* \otimes R_\rho\Big)^G,
  \end{gather*}
  see \cite[p.~267]{Craw2004}.  Here $(-)^D$ denotes the derived
  dual.  In particular,
  \begin{gather*}
    \Phi_\theta^{-1}(\sO \otimes R_\rho) = \sR_\rho.
  \end{gather*}
  Therefore, according to Theorem~\ref{thm:bkr},
  \begin{gather*}
    H^k(\sM_\theta,\sR_\rho^*\otimes\sR_\sigma)
     =H^k\big(\Hom_{D(\sM_{\theta})}(\sR_\rho,\sR_\sigma)\big) \\
  \hphantom{H^k(\sM_\theta,\sR_\rho^*\otimes\sR_\sigma)}{}
     \iso H^k\big(\Hom_{D^G(\C^3)}(\sO \otimes R_\rho,\sO\otimes R_\sigma)\big)
     = H^k\big(\C^3,\sO\otimes R_\rho^* \otimes R_\sigma\big)^G.
  \end{gather*}

  To prove the second part, we show the commutativity of the diagram
  \begin{equation*}
    \begin{tikzpicture}[baseline=(current  bounding  box.center)]
      \matrix (m) [matrix of math nodes,row sep=3em,column
      sep=3em,minimum width=2em ] {
        D(\sM_\theta) &
        D\big(\sM_\theta{\setminus}{\pi_{\theta}^{-1}(0)}\big)
        \\
        D^G\big(\C^3\big) &
        D^G\big(\C^3{\setminus}\{0\}\big). \\
      };
      \path[-stealth]
      (m-1-1) edge node [above] {$i^*$} (m-1-2)
      (m-2-1) edge node [below] {$j^*$} (m-2-2)
      (m-1-1) edge node [left] {$\Phi_\theta$} (m-2-1)
      (m-1-2) edge node [right] {$\Phi_{\sO_\Gamma}$} (m-2-2);
    \end{tikzpicture}
  \end{equation*}
  Here $\Phi_{\sO_\Gamma}$ is the Fourier--Mukai transform with kernel
  $\sO_\Gamma$, the structure sheaf of the graph of $\pi_{\theta}\colon
  \sM_\theta{\setminus}\pi_{\theta}^{-1}(0)\to \C^3{\setminus}\{0\}$.
  Note that under the identif\/ication
  $D^G(\C^3{\setminus}\{0\})=D((\C^3{\setminus}\{0\})/G)$,
  $\Phi_{\sO_\Gamma}$ becomes $(\pi_{\theta})_*$.

  Denote by $r$ and $s$ the projections from
  $\sM_{\theta}{\setminus}\pi_{\theta}^{-1}(0) \times
  \C^3{\setminus}\{0\}$ to $\sM_{\theta} {\setminus}\pi_{\theta}^{-1}(0)$
  and $\C^3{\setminus}\{0\}$, respectively.  Let $t\colon\sM_\theta\times
  \C^3{\setminus}\{0\}\to \C^3{\setminus}\{0\}$ denote the projection onto
  the second factor.  The following diagram summarises the situation:
 \begin{equation*}
    \begin{tikzpicture}[baseline=(current  bounding  box.center)]
      \matrix (m) [matrix of math nodes,row sep=3em,column
      sep=3em,minimum width=2em ] {
        & \sM_\theta{\setminus}\pi_{\theta}^{-1}(0)\times\C^3{\setminus}\{0\}
        \\
        \sM_\theta{\setminus}\pi_{\theta}^{-1}(0) &
        \sM_\theta\times\C^3{\setminus}\{0\}
        & \C^3{\setminus}\{0\}
        \\
        \sM_\theta & \sM_\theta\times\C^3 & \C^3 \\
      };
      \path[-stealth]
      (m-1-2) edge node [above] {$r$} (m-2-1)
              edge node [right] {$i\times\id$} (m-2-2)
              edge node [above] {$s$} (m-2-3)
      (m-2-1) edge node [left] {$i$} (m-3-1)
      (m-2-2) edge node [above] {$t$} (m-2-3)
              edge node [right] {$\id\times j$} (m-3-2)
      (m-2-3) edge node [right] {$j$} (m-3-3)
      (m-3-2) edge node [below] {$p$} (m-3-1)
              edge node [below] {$q$} (m-3-3);
    \end{tikzpicture}
  \end{equation*}
  It follows, essentially from the def\/inition of $\pi_{\theta}$, that
  \begin{gather}\label{eq:uog}
    (\id_{\sM_{\theta}}\times j)^* \sU_{\theta}
    = \big(i\times\id_{\C^3{\setminus} \{0\}}\big)_*\sO_\Gamma.
  \end{gather}
  Using \eqref{eq:uog} as well as the push-pull formula $(\sF\otimes
  f_*\sG) = f_*(f^*\sF\otimes \sG)$ we obtain
  \begin{gather*}
    j^*\circ\Phi_\theta(-)
     = j^*\circ q_*(p^*(-\otimes \rho_0)\otimes\sU_\theta)
    = t_* \big((\id_{\sM_{\theta}}\times j)^* p^*(-\otimes\rho_0)\otimes (\id_{\sM_{\theta}}\times j)^*\sU_\theta\big) \\
  \hphantom{j^*\circ\Phi_\theta(-)}{}
    = t_*\big((\id_{\sM_{\theta}}\times j)^* p^*(-\otimes\rho_0)\otimes (i\times\id_{\C^3{\setminus}\{0\}})_*\sO_\Gamma\big) \\
  \hphantom{j^*\circ\Phi_\theta(-)}{}
= t_*\big(i\times\id_{\C^3{\setminus} \{0\}}\big)_* \big((i\times j)^* p^*(-\otimes\rho_0)\otimes \sO_\Gamma\big) \\
  \hphantom{j^*\circ\Phi_\theta(-)}{}
 = s_* \big(r^*i^*(-\otimes \rho_0)\otimes \sO_\Gamma\big)
     = \Phi_{\sO_\Gamma}\circ i^*(-).
  \end{gather*}
  This concludes the proof.
\end{proof}

Before we embark on the proof of~\eqref{eq:h10}, it is useful to
recall some basic properties of local cohomology, see, e.g.,
\cite[Chapter~III, Exercise~2.3]{Hartshorne1977}.  Let $D$ be a~closed subset of~$X$ and let~$\cE$ be a sheaf on~$X$.  Denote by
$\Gamma_D(X,\sE)$ the subspace of $\Gamma(X,\sE)$ consisting of
sections whose support is contained in~$D$.  The functor
$\Gamma_D(X,-)$ is left-exact and its right derived functor
$H^\bullet_D(X,\sE):=R^\bullet\Gamma_D(X,\sE)$ is called the
\emph{local cohomology} of~$\cE$ with respect to~$D$.  Local
cohomology is related to the usual cohomology of $\cE$ by the
following long exact sequence
\begin{gather}\label{eq:les}
  \cdots \to H^i_D(X,\sE) \to H^i(X,\sE) \to H^i(X{\setminus}
  D,\sE|_{X{\setminus} D}) \stackrel{\delta}{\to} H^{i+1}_D(X,\sE) \to
  \cdots
\end{gather}
Moreover, it satisf\/ies \emph{excision}, that is, if $U$ is an open
subset in $X$ containing $D$, then there is a natural isomorphism
\begin{gather*}
  H^\bullet_D(X,\cE) \iso H^\bullet_D(U,\cE|_U).
\end{gather*}

\begin{proof}[Proof of Proposition~\ref{prop:av}]
  We have already reduced the proof of this to the proof of the
  vanishing statement~\eqref{eq:h10}.  Since by Proposition~\ref{prop:ci} we have $H^1(\sM_\theta,\sEnd(\sR))=0$, the long exact
  sequence associated to local cohomology yields
  \begin{gather*}
    H^0(\sM_\theta,\sEnd(\sR)) \stackrel{\delta}{\to}
    H_D^1\big(\bar\fM_\theta,\sEnd(\bar\sR)(-D)\big) \to
    H^1\big(\bar\fM_\theta,\sEnd(\bar\sR)(-D)\big) \to 0.
  \end{gather*}
  We show that the f\/irst map in this sequence is an isomorphism.  This
  gives the desired vanishing:
  $H^1(\bar\fM_\theta,\sEnd (\bar\sR )(-D))=0$.

  Let $H$ denote the hyperplane section in $\P^3$. By excision, we have
  \begin{gather*}
    H_D^1\big(\bar\fM_\theta,\sEnd (\bar\sR )(-D)\big)
    \iso H_D^1\big(\bar\fM_\theta{\setminus}
   \pi^{-1}(0),\sEnd (\bar\sR )(-D)\big) \\
   \hphantom{H_D^1\big(\bar\fM_\theta,\sEnd (\bar\sR )(-D)\big)}{}
    = H_H^1\big([\P^3{\setminus}\{[0:0:0:1]\}/G],\sO(-1)\otimes\End(R)\big) \\
   \hphantom{H_D^1\big(\bar\fM_\theta,\sEnd (\bar\sR )(-D)\big)}{}
    \iso H_H^1\big([\P^3/G],\sO(-1)\otimes\End(R)\big).
  \end{gather*}
  Here and in the following we omit to make the appropriate
  restriction of sheaves explicit, because confusion is unlikely to
  arise.  Using the above, we have the commutative diagram
  \begin{equation*}
    \begin{tikzpicture}[baseline=(current  bounding  box.center)]
      \matrix (m) [matrix of math nodes,row sep=3em,column
      sep=3em,minimum width=2em ] {
        H^0\big(\sM_\theta{\setminus}{\pi^{-1}(0)},\sEnd(\sR)\big) &
        H_D^1\big(\bar\fM_\theta,\sEnd(\sR)(-D)\big)
        \\
        H^0\big(\C^3{\setminus}\{0\}/G,\sO_{\C^3}\otimes\End(R)\big) &
        H_H^1\big([\P^3/G],\sO(-1)\otimes\End(R)\big). \\
      };
      \path[-stealth]
      (m-1-1) edge node [above] {$\tilde\delta$} (m-1-2)
      (m-2-1) edge (m-2-2)
      (m-1-1) edge node [left] {$\pi_*$} (m-2-1)
      (m-1-2) edge node [right] {$\iso$} (m-2-2);
    \end{tikzpicture}
  \end{equation*}
  We compose on the left with the commutative diagram in Proposition
  \ref{prop:ci}. Since $\delta = \tilde \delta \circ i^*$, we obtain
  the commutative diagram
  \begin{equation*}
    \begin{tikzpicture}[baseline=(current  bounding  box.center)]
      \matrix (m) [matrix of math nodes,row sep=3em,column
      sep=3em,minimum width=2em ] {
        H^0(\sM_\theta,\sEnd(\sR)) &
        H_D^1\big(\bar\fM_\theta,\sEnd (\bar\sR)(-D)\big) \\
        H^0\big(\C^3,\sO_{\C^3}\otimes\End(R)\big)^G &
        H_H^1\big([\P^3/G],\sO(-1)\otimes\End(R)\big). \\
      };
      \path[-stealth]
      (m-1-1) edge node [above] {$\delta$} (m-1-2)
      (m-2-1) edge (m-2-2)
      (m-1-1) edge node [left] {$\Phi_\theta$} (m-2-1)
      (m-1-2) edge node [right] {$\iso$} (m-2-2);
    \end{tikzpicture}
  \end{equation*}
  All the vertical arrows are isomorphisms. Moreover, by using the
  long exact sequence~\eqref{eq:les} and
  \begin{gather*}
    H^i\big([\P^3/G],\sO(-1)\otimes\End(R)\big)
      = H^i\big(\P^3,\sO(-1)\otimes\End(R)\big)^G = 0
  \end{gather*}
  for $i=0$ and $1$, it follows that the bottom map in the above
  diagram is an isomorphism.  Therefore, the map~$\delta$ must also be
  an isomorphism.
\end{proof}

\section[Dirac operators on $M_\theta$]{Dirac operators on $\boldsymbol{M_\theta}$}
\label{sec:dirac}

Let $(X,g)$ be an ALE spin manifold asymptotic to $\C^n/G$ and let $E$
be a complex vector bundle over~$X$ together with an asymptotically
f\/lat connection~$A$.  Denote by~$S^{\pm}$ be the spinor bundles on~$X$
and by $D^{\pm}_E$ the corresponding twisted Dirac operators.  We
denote by $W^{k,2}_{\delta}(X, S^\pm\otimes E)$ the completions of the
spaces of compactly supported sections with respect the weighted
Sobolev norm def\/ined by~\eqref{eq:wsf} using the covariant derivative
$\nabla_A$. Let $D^{\pm}_{E,\delta}\colon W^{k+1,2}_{\delta}(X,
S^{\pm}\otimes E) \to W^{k,2}_{\delta-1}(Z, S^{\mp} \otimes E)$ denote
the corresponding completion of the Dirac operator~$D^\pm_{E}$.

\begin{Theorem}
  For $\delta\in(-2n-1,0)$ the Dirac operator $D^\pm_{E,\delta}$ is
  Fredholm and its index is given by
  \begin{gather}\label{eq:index-ale}
    \ind D^{+}_{E,\delta} = \int_X \ch(E) \hat{A}(X) - \frac{\eta_E}{2}.
  \end{gather}
  Here  $\ch(E)$ denotes the Chern character of $E$ as a differential form,
  $\hat{A}(X)$ denotes the Hirzebruch $\hat{A}$-polynomial
  applied to the Pontryagin forms $p_i(X)$ of the ALE metric on $X$,
   $\eta_E(s):=\sum\limits_{\lambda\neq 0}\sign(\lambda) |\lambda|^{-s}$ is
  the $\eta$-function of the spectrum of the Dirac operator restricted to
  the boundary at infinity $S^{2n-1}/G$ of the ALE manifold $X$, and
  $\eta_E:=\eta_E(0)$ is the $\eta$-invariant.
\end{Theorem}

\begin{proof}
  The fact that $D^\pm_{E,\delta}$ is Fredholm is proved as in
  Proposition~\ref{prop:Laplace} by noting that the set of indicial
  roots does not intersect~$(-2n-1, 0)$.  This can be seen, for
  example, by realizing that the indicial roots correspond to the
  eigenvalues of the Dirac operator on $S^{2n-1}/G$ shifted by~$-\frac{2n-1}{2}$.  The index formula follows from
  Atiyah--Patodi--Singer index theorem~\cite{Atiyah1975a}.
\end{proof}

If $(X,g)$ is a K\"ahler manifold, then there is a one-to-one
correspondence between spin structures on $X$ and holomorphic square
roots of the canonical line bundle~$K_X$,
see \cite[Theorem~2.2]{Hitchin1974}. We are interested in the case
when $X = M_{\theta}$ for $\theta \in \Theta_{\Q}$ a generic stability
parameter.  Since $M_\theta$ is a crepant resolution of~$\C^3/G$, its
canonical line bundle is holomorphically trivial.  In par\-ti\-cular, for
any K\"ahler metric on $M_\theta$ there is a canonical spin structure
corresponding to taking the square root of the canonical bundle to be
the trivial holomorphic line bundle on~$M_\theta$. The corresponding
spinor bundles are
\begin{gather*}
  S^{+} = \Lambda^{0,\text{even}} T_\C^*M_\theta \qquad\text{and}\qquad
  S^{-} = \Lambda^{0, \text{odd}} T_\C^*M_\theta.
\end{gather*}
We f\/ix this spin structure on $M_\theta$ for the rest of the section.
Suppose that $A$ is a $\U(n)$-connection on a holomorphic bundle~$\cE$
compatible with the given holomorphic structure.  If the metric on~$M_\theta$ induces the product connection on our chosen square root of~$K_{M_\theta}$,
i.e., the metric is Ricci-f\/lat, then the corresponding twisted Dirac
operator is given by
\begin{gather}\label{eq:D-1}
  D^\pm_{\cE} = \sqrt2 \big(\delbar_A + \delbar_A^*\big).
\end{gather}

We now show that for special choices of $\cE$ involving the
tautological bundles on $M_\theta$, the index of the corresponding
weighted Dirac operator is zero.  Recall that by Theorem~\ref{Thm_KRF_HYM},
there exists an ALE Calabi--Yau metric $g_{\theta, \RF}$ on $M_\theta$
in the K\"ahler class of~$g_{\theta}$ and a HYM connection on each of
the tautological holomorphic line bundles $\cR_{\rho}$ with $\rho\in
\Irr(G)$.  Moreover these induce a~HYM connection $A$ on $\cR =
\bigoplus_{\rho\in \Irr(G)} \cR_{\rho}$ which is inf\/initesimally
rigid.

\begin{Proposition}\label{prop:mv}
  Let $\theta\in\Theta_\Q$ be generic.
  \begin{enumerate}\itemsep=0pt
  \item[$(1)$] If $g$ is an ALE K\"ahler metric on $M_\theta$ in the same
    K\"ahler class as $g_\theta$, then
    \begin{gather*}
      \ind D^{+}_{\cR_\rho\otimes \cR_\sigma^*, \delta}=0,
    \end{gather*}
    for all $\rho, \sigma \in \Irr(G)$ and for all $\delta \in (-5,
    0)$.
  \item[$(2)$] If $g=g_{\theta,\RF}$ and $\cR$ is equipped with the HYM
    connection $A$ given by Theorem~{\rm \ref{Thm_KRF_HYM}(3)}, then
    the twisted Dirac operator $D^{+}_{\cEnd(\cR),\delta}$ is an
    isomorphism for all $\delta \in (-5, 0)$.
  \end{enumerate}
\end{Proposition}

\begin{proof}
  The index of $D^{+}_{\cR_{\rho}\otimes \cR_{\sigma}^*, \delta}$ is
  unchanged under deforming the ALE K\"ahler metric and the
  asymptotically f\/lat connection on $\cR_\rho\otimes\cR_\sigma^*$.
  Hence there is no loss in assuming that $g=g_{\theta,\RF}$ and that
  $\cR_\rho\otimes\cR_\sigma^*$ has been equipped with its HYM metric.
  Since the twisted Dirac operators preserve the holomorphic splitting
  $\cEnd(\cR) = \bigoplus_{\rho,\sigma\in \Irr(G)}\cR_{\rho}^*\otimes
  \cR_{\sigma}$, we are thus left with proving the second statement.

  First we prove that $\ind D^+_{\cEnd(\cR),\delta}=0$.  Let $\Omega$
  be a nowhere vanishing holomorphic volume form on $M_\theta$ and let
  $*\colon \Lambda^{p,q} T^*M_\theta \otimes \cE \to \Lambda^{3-p,3-q}
  T^*M_\theta \otimes \cE^*$ denote the Hodge-$*$-operator for some
  holomorphic bundle $\cE$.  Then we have an isomorphism of vector
  bundles $S^{+} \otimes \cE \iso S^{-} \otimes \cE^*$, given by
  \begin{gather*}
    S^+\otimes \cE=\Lambda^{0,\text{even}} T^*M_\theta\otimes \cE
    \stackrel{\Omega\wedge -}{\longrightarrow} \Lambda^{3,{\rm even}}
    T^*M_\theta \otimes \cE \\
\hphantom{S^+\otimes \cE}{}
    \stackrel{*}{\longrightarrow}
    \Lambda^{0,{\rm odd}} T^*M_\theta \otimes \cE^* = S^-\otimes
    \cE^*.
  \end{gather*}
  Similarly, $S^-\otimes \cE \iso S^+\otimes \cE^*$.  These
  isomorphisms identify $D^+_\cE$ with $D^-_{\cE^*}$.  Consequently
  \begin{gather*}
    \ind D^+_{\cE,\delta}
    =\ind D^-_{\cE^*,\delta}.
  \end{gather*}
  Moreover, the $L^2$-adjoint of $D^{-}_{\cE^*, \delta}$ is
  $D^{+}_{\cE^*,-5-\delta}$ and thus
 \begin{gather*}
   \ind D^{-}_{\cE^*,\delta} = - \ind D^{+}_{\cE^*, -5-\delta}.
 \end{gather*}
 For $\cE = \cEnd(\cR)$ and $\delta=-\frac52$, the two identities
 above give that $\ind D^+_{\cEnd(\cR),-\frac{5}{2}}=0$.  Since the
 index is constant for $\delta \in (-5, 0)$, we must have $\ind
 D^+_{\cEnd(\cR),\delta}=0$ for all $\delta \in (-5, 0)$.

 To complete the proof, we show that $\coker
 D^+_{\cEnd(\cR),\delta}=0$ for $\delta\in(-5,0)$ or, equivalently,
 that $\ker D^{-}_{\cEnd(\cR),-5-\delta}=0$.  By abuse of notation we
 also denote by $A$ the connection induced by~$A$ on~$\cEnd(\cR)$.
 Each $(\phi_1, \phi_3) \in \ker D^{-}_{\cEnd(\cR)}$ satisf\/ies
 \begin{gather}\label{eq:mv-5}
   \delbar_A^* \phi_1 =0 \qquad\text{and}\qquad \delbar_A \phi_1 +
   \delbar_A^*\phi_3 =0.
 \end{gather}
 The second identity gives
 \begin{gather}\label{eq:mv-2}
   \delbar_A \delbar_A^* \phi_3 =0.
 \end{gather}
 Arguing as in Proposition~\ref{prop:decay} one shows that $\phi_3 =
 O(r^{-4})$ and $\delbar^*_A \phi_3 = O(r^{-5})$.  Hence, taking
 the $L^2$-inner-product with $\phi_3$ in~\eqref{eq:mv-2} and
 integrating by parts, we conclude that $\delbar_A^* \phi_3 =0$.  On
 the other hand, we have the splitting $\nabla_A = \partial_A +
 \delbar_A$ and the K\"ahler identity gives $\partial_A = i [\Lambda,
 \delbar^*_A]$. Using this, the above yields $\partial_A \phi_3 =0$.
 If we write $\phi_3 = \bar\Omega\otimes s$ with $s$ a smooth section
 of $\cEnd(\cR)$, then $\del_{A} (\bar\Omega \otimes s) = \bar\Omega
 \wedge \del_{A}s =0$ and thus $\del_{A}s=0$.  Since the connection
 $A$ on $\cR$ is HYM and $g_{\theta,\RF}$ is Ricci-f\/lat, the
 Weitzenb\"ock formula gives $2\Delta_{\del_A}=\nabla_A^*\nabla_A$.  It
 follows that $s$ is parallel and hence must vanish, as it is zero at
 inf\/inity.  This implies that $\phi_3 =0$. From~\eqref{eq:mv-5} we
 deduce that $\phi_1$ satisf\/ies $\delbar_A\phi_1 = 0$ and $\delbar_A^*
 \phi_1 =0$.  Then, by Lemma~\ref{lem:vanishing}, $\phi_1$ must vanish
 identically.
\end{proof}

\section{The proof of Theorem~\ref{Thm_MultiplicativeIdentity}}

\label{sec:mckay}
\label{Sec_AnalyticProof}

We now prove Theorem~\ref{Thm_MultiplicativeIdentity}.  We start by discussing the
sense in which formula~\eqref{Eq_MultiplicativeIdentity} is valid. Since
$M_{\theta}$ is non-compact, this formula cannot be interpreted as a
dif\/ference of triple products in $H^*(M_{\theta},\R)$ as
$H^6(M_{\theta},\R) =0$.

However, each of the tautological line bundles $\cR_{\rho}$
have natural asymptotically f\/lat connec\-tions~$A_{\theta, \rho}$
compatible with their holomorphic structure. Via Chern--Weil theory,
$c_1(\cR_{\rho})$ is represented in $H^2(M_{\theta},\R)$ by the de
Rham cohomology class of $ \frac{i}{2\pi}F_{A_{\theta,\rho}}$. These
are the dif\/ferential forms used in the Atiyah--Patodi--Singer
formula~\eqref{eq:index-ale} to represent the Chern character of the
bundle~$\cR_{\rho}$.  Using this formula we will show below how to
obtain~\eqref{Eq_MultiplicativeIdentity} (in fact, the equivalent formulation~\eqref{Eq_MultiplicativeIdentity2}) at the level of dif\/ferential forms
representing the f\/irst Chern classes.

Moreover, since we are on an ALE manifold, the long exact sequence in
cohomology
\begin{gather*}
  \cdots \to H^1\big(S^5/G,\R\big) \to H^2_{c}(M_{\theta},\R) \stackrel{j}\to
  H^2(M_{\theta},\R) \to H^2\big(S^5/G,\R\big) \to \cdots
\end{gather*}
gives that the homomorphism $j$, which takes the class of a compactly
supported form in $H^2_c(M_{\theta},\R)$ to its de Rham
representative in $H^2(M_{\theta},\R)$, is an isomorphism. As such,
there exists a compactly supported $2$-form $\alpha_{\rho}$ so that
\begin{gather*}
  \tfrac{i}{2\pi}F_{A_{\theta,\rho}} = \alpha_{\rho} + d\beta_{\rho}
\end{gather*}
with $\beta_{\rho}$ a $1$-form on $M_{\theta}$. Since
$F_{A_{\theta,\rho}} = dA_{\theta, \rho}$, after possibly modifying
$\alpha_{\rho}$, we can take $\beta_{\rho}$ to be $\frac{i}{2\pi}
A_{\theta, \rho}$ on the ALE end of~$M_{\theta}$. From
Theorem~\ref{thm:asympt}(3) we know that the
connections $A_{\theta,\rho}$ are asymptotically f\/lat of order $1$ and
that the curvature $F_{A_{\rho}} = O(r^{-4})$ as $r\to \infty$. Then
Stokes' formula gives
\begin{gather*}
  \int_{M_{\theta}} c_1(\cR_{\rho}) c_1(\cR_{\sigma}) c_1 (\cR_{\tau})
  = \int_{M_{\theta}} \alpha_{\rho} \wedge \alpha_{\sigma} \wedge
  \alpha_{\tau}
\end{gather*}
for all $\rho, \sigma, \tau \in \Irr(G)$.  From here, we see
that~\eqref{Eq_MultiplicativeIdentity} is interpreted topologically in terms of
triple products
\begin{gather*}
  \int_{M_\theta}\colon \ H^2_c(M_{\theta},\R) \times H^2_c(M_{\theta},\R) \times H^2_c(M_{\theta},\R) \to \R,
\end{gather*}
with the f\/irst Chern classes $c_1(\cR_{\rho})$ thought of as their own
images in $H^2_c(M_{\theta},\R)$ under $j^{-1}$.

Therefore, it remains to prove Theorem~\ref{Thm_MultiplicativeIdentity} at the level of
dif\/ferential forms.  As we already mentioned, the proof uses the
Atiyah--Patodi--Singer index theorem for ALE
manifolds~\eqref{eq:index-ale}. In order to apply it, we need to
compute the $\eta$-invariant term that appears in this formula.

\begin{Proposition}
  Let $G$ be a finite subgroup of $\SL(n,\C)$ acting freely on $\C^n
  {\setminus} \{0\}$. Assume that~$X$ is a smooth ALE spin manifold
  asymptotic to~$\C^n/G$ and let~$(E,A)$ be a asymptotically flat
  bundle on~$X$ whose fiber at infinity is~$E_{\infty}$. Then, the
  $\eta$-invariant for the Dirac operator~$D_{E,\delta}$ on~$X$ is
  given by
  \begin{gather}\label{eq:eta}
    \eta_E = -\frac{2}{|G|} \sum_{\substack{g \in G \\ g \neq I_n}}
    \frac{\chi_{E_{\infty}} (g)}{\sum\limits_{i=0}^n (-1)^i \chi_{\Lambda^i\C^n}(g)},
  \end{gather}
  provided $-2n+1<\delta <0$.  In this formula $\chi_{E_{\infty}}$ denotes the
  character of the representation corresponding to the action of~$G$
  on the vector space~$E_{\infty}$.
\end{Proposition}

\begin{Remark}
  Note that for any $g\in \SL(n,\C)$, $\sum\limits_{i=0}^n (-1)^i
  \chi_{\Lambda^i\C^n} (g) = \det(\id_{\C^n}-g)$. Since $G$ is chosen to act
  freely on $\C^n$, $\det(I_n-g) \neq 0$ for all $g \in G{\setminus}
  \{I_n\}$, and thus all the denominators in formula~\eqref{eq:eta} are
  non-zero.
\end{Remark}

This proposition is a consequence of the Lefschetz f\/ixed-point
formula, in the sense that~$\eta_E$ is the contribution from the f\/ixed
locus under the action of~$G$ on~$\C^n$. It can be also proved using
the def\/inition of the $\eta$-invariant as the analytic continuation at~$0$ of the $\eta$-series correspon\-ding to the spectrum of the Dirac
operator on the boundary at inf\/inity of the orbi\-fold~$\C^n/G$.  This
last approach gives the generalisation of the above formula to the
case of non-isolated singularities~\cite{Degeratu2001}.

Let $G\subset\SL(n,\C)$ be a f\/inite subgroup.  Then for each $\rho\in
\Irr(G)$ we have the decomposition into irreducibles
\begin{gather}\label{eq:vr}
  \Lambda^i\C^n \otimes \rho
    = \sum_{\sigma \in \Irr(G)} a^{(i)}_{\rho\sigma}\ \sigma
\end{gather}
with $a^{(i)}_{\rho\sigma}\in\N_0$.  As in Ito and
Nakajima~\cite{Ito2000}, we def\/ine
\begin{gather}\label{eq:cartan}
  c_{\rho\sigma} := \sum_{i=0}^n (-1)^i a_{\rho\sigma}^{(i)}
\end{gather}
and set
\begin{gather}
  \label{Eq_C}
  \tilde C:=(c_{\rho\sigma})_{\rho,\sigma \in \Irr(G)} \qquad\text{and}\qquad
  C:=(c_{\rho\sigma})_{\rho,\sigma \in \Irr_0(G)}.
\end{gather}

\begin{Remark}
  When $G$ is a f\/inite subgroup of $\SL(2,\C)$, the matrix~$C$ is the
  Cartan matrix of the simple Lie algebra corresponding to~$G$, while~$\tilde{C}$ is its extended version.  This is the essence of the classical
  McKay correspondence~\cite{McKay1980}.  For $n\geq 3$ this matrix has, in general, entries on the diagonal which are not equal to~$2$, and hence is
  neither the Cartan matrix associated to a~Lie algebra nor is it a~generalised Cartan matrix as appearing in the context of Lie
  algebras.
\end{Remark}

If $n$ is even, then the virtual representation $\sum\limits_{i=0}^n(-1)^i
\Lambda^i \C^n$ is self-dual and $C$ is symmetric; otherwise the
virtual representation is anti-self-dual and $C$ is anti-symmetric.
If~$G$ is abelian, as is the case in our situation, then every
irreducible representation is one-dimensional and thus $\sum\limits_{\sigma
  \in \Irr(G)} a_{\rho\sigma}^{(i)} = \dim  \Lambda^i\C^n$.  Combined
with formula~\eqref{eq:cartan} this gives
\begin{gather}\label{eq:c}
  \sum_{\sigma \in \Irr(G)} c_{\rho\sigma} = 0
\end{gather}
for all $\rho \in \Irr(G)$.

\begin{proof}[Proof of Theorem~\ref{Thm_MultiplicativeIdentity}]
  By Proposition~\ref{prop:mv} we know that $\ind
  D^{+}_{\cR_{\rho}\otimes \cR_{\sigma}^*, \delta} =0$ for all weights
  $\delta \in (-5, 0)$ and for all $\rho, \sigma \in \Irr(G)$.  Thus
  the index formula~\eqref{eq:index-ale} takes the form
  \begin{gather*}
    \int_{M_\theta} \ch(\cR_{\rho} \otimes \cR_{\sigma}^*) \hat
    A(M_\theta) = \frac{\eta_{R_{\rho} \otimes R_{\sigma^*}}}{2}.
  \end{gather*}
  Multiplying this equation by the matrix $\widetilde{C}$ yields
  \begin{gather}\label{eq:mckay-2-1}
    \sum_{\rho \in \Irr(G)} c_{\tau \rho} \int_{M_\theta}
    \ch(\cR_{\rho} \otimes \cR_{\sigma}^*) \hat A(M_\theta) =
    \sum_{\rho \in \Irr(G)} c_{\tau \rho} \frac{\eta_{R_{\rho} \otimes
        R_{\sigma^*}}}{2}
  \end{gather}
  for all $\tau \in \Irr(G)$.  The left-hand side
  of~\eqref{eq:mckay-2-1} can be written as
  \begin{gather*}
    \sum_{\rho\in\Irr(G)} c_{\tau\rho}\int_{M_\theta} \tilde\ch(\cR_{\rho})
         \tilde\ch(\cR_{\sigma}^*) \\
  \qquad\quad{}  +\sum_{\rho \in \Irr(G)} c_{\tau\rho}\int_{M_\theta}\ch(\cR_{\rho}) \hat{A}(M_\theta)
    + \sum_{\rho\in \Irr(G)} c_{\tau\rho}\int_{M_\theta}\ch(\cR_{\sigma}^*)\hat A (M_\theta)
 \\
 \qquad{}   = \sum_{\rho\in\Irr(G)} c_{\tau\rho}\int_{M_\theta} \tilde\ch(\cR_{\rho})
         \tilde\ch(\cR_{\sigma}^*)
    + \sum_{\rho\in \Irr(G)}c_{\tau\rho} \frac{\eta_{R_{\rho}}}{2}
    + \sum_{\rho\in \Irr(G)} c_{\tau\rho} \frac{\eta_{R_{\sigma}^*}}{2}
  \end{gather*}
  for all $\tau \in \Irr(G)$.  Here we have used the fact that $\hat A
  (M_\theta) = 1 + \hat A_4 (M_\theta)$ since $M_\theta$ is
  $6$-dimensional.  By~\eqref{eq:c} the third term in the above
  expression vanishes.  Therefore we can rewrite~\eqref{eq:mckay-2-1}
  as
  \begin{gather}\label{eq:mckay-1}
    \sum_{\rho\in\Irr(G)} c_{\tau\rho}\int_{M_\theta}
    \tilde\ch(\cR_{\rho}) \tilde\ch(\cR_{\sigma}^*) =
  \sum_{\rho \in\Irr(G)} c_{\tau\rho}  \frac{\eta_{R_{\rho} \otimes R_{\sigma^*}}}{2} - \sum_{\rho\in
      \Irr(G)}c_{\tau\rho} \frac{\eta_{R_{\rho}}}{2} .
  \end{gather}
  Taking the characters of \eqref{eq:vr} and summing over $i$ with
  alternating signs we obtain
  \begin{gather*}
    \big(\chi_{\C^3} (g) - \chi_{\Lambda^2 \C^3} (g)\big) \chi_{\tau} (g) = -
    \sum_{\rho \in \Irr(G)} c_{\tau\rho} \chi_{\rho} (g),
  \end{gather*}
  which gives
  \begin{gather*}
    \chi_{\tau\otimes \sigma^*} (g) = - \sum_{\rho\in \Irr(G)} c_{\tau
      \rho} \frac{\chi_{\rho\otimes\sigma^*}(g)}{\chi_{\C^3}(g) -
      \chi_{\Lambda^2 \C^3}(g)}
  \end{gather*}
  for all $g \in G{\setminus} \{e\}$.  Summing over all such $g$ and
  using formula \eqref{eq:eta} for the $\eta$-invariant, we obtain
  \begin{gather*}
    - 2 \left(\delta_{\tau\sigma}- \frac{1}{|G|}\right) = \sum_{\rho \in
      \Irr(G)} c_{\tau\rho} \eta_{\cR_{\rho} \otimes \cR_{\sigma}^*}
  \end{gather*}
  Hence, \eqref{eq:mckay-1} yields
  \begin{gather*}
    \sum_{\rho\in\Irr(G)} c_{\tau\rho}\int_{M_\theta}
    \tilde\ch(\cR_{\rho}) \tilde\ch(\cR_{\sigma}^*) =
    -\left(\delta_{\tau\sigma}-\frac{1}{|G|}\right) + \left(\delta_{\tau \rho_0}
    -\frac{1}{|G|}\right)
  \end{gather*}
  for all $\tau, \sigma \in \Irr(G)$.  For $\tau \in \Irr_0(G)$ this
  becomes
  \begin{gather}\label{eq:cartan_invertible}
    \sum_{\rho\in\Irr(G)} c_{\tau\rho}\int_{M_\theta}
    \tilde\ch(\cR_{\rho}) \tilde\ch(\cR_{\sigma}^*) =
    -\delta_{\tau\sigma}.
  \end{gather}
  Since $\tilde\ch(\cR_{\rho_0})=0$, it follows that the matrix $C =
  (c_{\tau\rho})_{\tau,\rho \in \Irr_0(G)}$ is invertible and
  \begin{gather}\label{eq:mult-3}
    \int_{M_\theta} \tilde\ch(\cR_{\rho}) \tilde\ch(\cR_{\sigma}^*) =
    -\big(C^{-1}\big)_{\rho\sigma},
  \end{gather}
  for all $\rho,\sigma\in \Irr_0(G)$, which is
  precisely~\eqref{Eq_MultiplicativeIdentity}.
\end{proof}

\begin{Remark}
  Note that formula~\eqref{eq:cartan_invertible} shows that the matrix~$C$ is invertible. In the case of a f\/inite subgroup of~$\SL(2,\C)$
  the invertibility of~$C$ was a direct consequence of the McKay
  correspondence, given that~$C$ is the Cartan matrix associated to a
  simply-laced Coxeter--Dynkin diagram~\cite{McKay1980}.
\end{Remark}

\begin{Remark}
  From formula~\eqref{eq:mult-3} it can be easily deduced that the set
  $\{\tilde{\ch}(\cR_{\rho})\colon \rho \in \Irr_0(G)\}$ is linearly
  independent in $H^*_c(M_{\theta},\R)$. To see this, note that if we
  had a linear combination $\sum\limits_{\rho \in \Irr_0(G)} a_{\rho}
  \tilde{\ch}(\cR_{\rho}) =0$, then by multiplying with
  $\tilde\ch(\cR_{\sigma}^*)$ for $\sigma \in \Irr_0(G)$ the
  left-hand side becomes $-\sum\limits_{\rho \in \Irr_0(G)} a_{\rho}
  (C^{-1})_{\rho\sigma}$ and thus all $a_{\rho}$ must vanish.
\end{Remark}

\subsection*{Acknowledgments}
A.D.~would like to thank Tom Mrowka, Tam\'as Hausel, Rafe Mazzeo and Mark Stern for useful conversations about dif\/ferent aspects of this work.
A.D.~was supported by the DFG via SFB/Transregio 71 ``Geometric Partial Dif\/ferential Equations''.
Parts of this article are the outcome of work undertaken by~T.W.\ while working on his PhD thesis at Imperial College London, supported by European Research Council Grant 247331.
T.W.~would like to thank his supervisor Simon Donaldson for his support.
Both authors would like to thank the anonymous referee of an earlier version of this article for pointing out a way of deriving the multiplicative formula~\eqref{Eq_MultiplicativeIdentity2} from the work of Ito and Nakajima~\cite{Ito2000}.

\pdfbookmark[1]{References}{ref}
\LastPageEnding

\end{document}